\documentclass[a4paper, 12pt]{article}

\usepackage{amssymb}
\usepackage{amsmath}
\usepackage{amsthm}
\usepackage{url}
\usepackage{mathrsfs}
\usepackage{graphicx}
\usepackage{sectsty}
\usepackage[letterpaper,hmargin=1.0in,vmargin=1.0in]{geometry}
\usepackage[colorlinks,linkcolor=black,citecolor=black,filecolor=black,urlcolor=black]{hyperref}
\usepackage{color}

\usepackage[british]{babel}

\usepackage{xspace}
\usepackage{calc}
\usepackage{euscript}
\usepackage{amsfonts}
\usepackage{enumerate}

\usepackage[textsize=tiny]{todonotes}

\usepackage{a4wide}
\parskip \medskipamount

\sectionfont{\large}
\subsectionfont{\normalsize}
\numberwithin{equation}{section}

\newtheorem{maintheorem}{Theorem}
\newtheorem{theorem}{Theorem}[section] 
\newtheorem{lemma}[theorem]{Lemma}
\newtheorem{propo}[theorem]{Proposition}
\newtheorem{fact}[theorem]{Fact}

\newtheorem{example}[theorem]{Example}
\newtheorem{corol}[theorem]{Corollary}

\newtheorem{remark}[theorem]{Remark}
 
  \newtheorem{problem}[theorem]{Problem}

\def\P{\mathbb{P}}

\newcommand{\Ind}[1]{\mathbf{1}_{  \{#1\} }}

\newcommand{ \hit}{ t_{\mathrm{hit}} }
\newcommand{\E}{{\mathbb{E}}}

\renewcommand{\Pr}{ \mathbb P}

\newcommand{ \rel}{ t_{\mathrm{rel}} }
\newcommand{ \mix}{ t_{\mathrm{mix}} }
\newcommand{ \mixin}{ t_{\mathrm{mix}}^{(\infty)} }
\newcommand{ \mixtv}{ t_{\mathrm{mix}}^{\mathrm{TV}} }

\newcommand{ \TV}{ \mathrm{TV} }

\newcommand{\eps}{\epsilon}

\newcommand{\la}{\lambda}

\DeclareMathSymbol{\leqslant}{\mathalpha}{AMSa}{"36} 
\DeclareMathSymbol{\geqslant}{\mathalpha}{AMSa}{"3E} 
\DeclareMathSymbol{\eset}{\mathalpha}{AMSb}{"3F}     
\renewcommand{\le}{\;\leqslant\;}                   
\renewcommand{\ge}{\;\geqslant\;}                   

\renewcommand{\epsilon}{\varepsilon}

\newcommand{\sfrac}[2]{\mbox{\small $\frac{#1}{#2}$}}
\newcommand{\ssfrac}[2]{\mbox{\footnotesize $\frac{#1}{#2}$}}
\newcommand{\half}{\ssfrac{1}{2}}

\newcommand{\N}{\mathbb N}
\newcommand{\R}{\mathbb R}

\begin{document}

\title{Some inequalities for reversible Markov chains and branching random walks via spectral optimization}
\author{ Jonathan Hermon
\thanks{
University of British Columbia, Vancouver, UK. E-mail: {\tt jhermon.math.ubc.ca}. Financial support by
the EPSRC grant EP/L018896/1 and by NSERC discovery and accelerator  grants.}}
\date{}
\maketitle

\begin{abstract}

We  present results relating mixing times to the intersection time of branching random walk (BRW) in which the logarithm of the expected number of particles grows at rate of the spectral-gap  $\mathrm{gap}$ . This is a finite state space analog of a critical branching process. Namely, we show that the maximal expected hitting time of a state by such a BRW is up to a universal constant larger than the $L_{\infty}$ mixing-time, whereas under transitivity the same is true for the intersection time of two independent such BRWs.    

Using the same methodology, we show that for a sequence of reversible Markov chains, the  $L_{\infty}$ mixing-times $t_{\mathrm{mix}}^{(\infty)} $ are of smaller order than the maximal hitting times $t_{\mathrm{hit}}$ iff the product of the spectral-gap and $t_{\mathrm{hit}}$ diverges, by establishing the inequality $t_{\mathrm{mix}}^{(\infty)} \le \frac{1}{\mathrm{gap}}\log(et_{\mathrm{hit}} \cdot \mathrm{gap}) $. This resolves a conjecture of Aldous and Fill \cite{aldousfill} Open Problem 14.12 asserting that under transitivity the condition that $ t_{\mathrm{hit}}  \gg \frac{1}{\mathrm{gap}} $ implies mean-field behavior for the coalescing time of coalescing random~walks.
\end{abstract}

\paragraph*{\bf Keywords:}
{\small Mixing times, hitting times, spectral optimization,   vertex-transitive graphs, intersection times, branching random walk, spectral-gap, coalescing random walk.
}


\section{Introduction}
The purpose of this article is to establish  new relations between certain natural quantities associated with a reversible Markov chain on a finite state space using a new `spectral optimization' technique. These quantities are the maximal  hitting time of a state $\hit$, the $L_{\infty}$ (a.k.a.\ uniform) mixing time $\mixin$ and the spectral-gap $\mathrm{gap}$. One of these relations (implied by Theorem \ref{p:mixequiv}) is the inequality $\mixin \le \frac{1}{\mathrm{gap}}\log(et_{\mathrm{hit}} \cdot \mathrm{gap})$, which as we explain below   is the missing component which previously remained open concerning a conjecture of Aldous and Fill about a sufficient condition for mean field behavior for coalescing random walk.
We then give in \S\ref{s:higherorder} (Theorem \ref{thm:higherorder}) higher order versions of this involving a parameter $\ell \in \N $, with $\ell=1$ corresponding to hitting times as in Theorem \ref{p:mixequiv} and $\ell=2$ to intersection times, as in Theorem \ref{thm:BRWintro2}. The connection of the cases $\ell=1,2$ with branching random walk which splits at a rate of order of the spectral-gap is given in \S\ref{s:interintro} in Theorems  \ref{thm:BRWintro} and \ref{thm:BRWintro2}. 
In particular, we give an interpretation to the quantity $\frac{1}{\mathrm{gap}}\log(et_{\mathrm{hit}} \cdot \mathrm{gap})$ (from the above bound on $\mixin$) in terms of hitting times for such a branching random walk, as well as a probabilistic interpretation in the transitive setup of the bound of Theorem \ref{thm:higherorder} in the case $\ell=2$  in terms of intersection times of two independent such branching random walks. 

The aforementioned relations refine previously known relations, such as the classical inequality \eqref{e:classicthm1} below as well as the work of Peres et al.\ \cite{inter} which was the first to relate intersection times to mixing times (this work refines \cite{inter} only in the transitive\footnote{We say a Markov chain on a countable state space $V$  with transition matrix $P$ is  \emph{transitive} if for every $x,y \in V$ there is a bijection $f:V \to V$ such that $f(x)=y$ and $P(x,z)=P(y,f(z))$ for all $z \in V$.
} setup).

It is a classical result of Aldous \cite{aldous1982some} that the total variation mixing time of a reversible Markov chain can be characterized in terms of hitting times of sets. For  recent refinements of this result see \cite{PS,Olivehit,basu}; see also \cite{L2} for an analogous result concerning the $L_{\infty}$ mixing time. It is also well-known that the $L_{\infty}$  mixing time can be (usually very loosely) bounded   in terms of the maximal (expected) hitting time of a state (see \eqref{e:mixlehit}) and also (usually less loosely)  in terms of the spectral-gap \eqref{e:trelmix}  (see \eqref{e:classicthm1} below for a combination of parts of \eqref{e:trelmix}-\eqref{e:mixlehit}). Our Theorem \ref{p:mixequiv} gives a more precise relation between these quantities than these two bounds. This new relation is also of theoretical interest, as it gives a simple spectral characterization for the condition that the mixing time and the maximal hitting time are of the same order. The implication of this for a conjecture of Aldous and Fill about mean field behavior of coalescing random walk is discussed later.

 Let $(X_t )_{t \ge 0}$ be an irreducible Markov chain on a finite state space $V$ with transition matrix $P$ and stationary distribution $\pi$ (we use this notation throughout the paper). We  consider the continuous-time rate 1 version of the chain. Let $H_t:=e^{-t(I-P)}$ be its heat kernel  (so that $H_t(\cdot,\cdot \cdot)$ are the time $t$ transition probabilities).
We note that our results are valid also in the discrete-time setup when $\min_{x\in V} P(x,x)$ is bounded away form 0. We  denote the eigenvalues of the Laplacian $I-P$ by $0=\la_1<\la_2 \le \ldots \la_{|V|} \le 2 $. The \emph{spectral-gap}   $\mathrm{gap} $  is defined as $\la_2$ and the \emph{relaxation-time} as $\rel:=\frac{1}{\la_2}$. 

 The  maximal expected \emph{hitting time} of a \textcolor{black}{state} is given by  \[\hit:=\max_{x,y \in V }\E_x[T_y], \qquad \text{where} \qquad T_y:=\inf\{t:X_t=y \}. \]  
It is well-known that the \emph{random target time}  $t_{\odot}:=\sum_{y} \pi(y)\E_x[T_y] $ is independent of $x$ (see \S\ref{s:hitbackground}; the notation  $t_{\odot}$  is borrowed from~\cite{{LPW}}) and so it is equal to the \emph{average hitting time}  \[t_{\odot}=\sum_{x,y}\pi(y)\pi(x) \E_{x}[T_y]=\sum_{y}\pi(y)\E_{\pi}[T_y]=\sum_{y}\pi(y)t_{\pi \to y}, \]
where $t_{\pi \to y}:=\E_{\pi}[T_y]$ is the expected hitting time of a state $y$ starting at equilibrium. 

We denote the $\varepsilon$ $L_{\infty}$ \emph{mixing time} and the $\varepsilon$  \emph{average} $L_2$ \emph{mixing time}, respectively, by\footnote{See \eqref{eq: taupeps}-\eqref{e:aveL2} and the discussion in \S\ref{s:mix} for the second equality for $\mix^{(\infty)}$ as well as for some motivation for $t_{\mathrm{ave-mix}}^{(2)}$  and an  interpretation of it as an ``average $L_2$ mixing time" (as suggested by its name).}  \[\mix^{(\infty)}(\varepsilon):=\inf\{t:\max_{a,b} |\sfrac{H_t(a,b)}{\pi(b)}-1| \le \varepsilon \}= \inf\{t:\max_{a} \sfrac{H_t(a,a)}{\pi(a)}-1 \le \varepsilon \},\]  \[t_{\mathrm{ave-mix}}^{(2)}(\varepsilon ):=\inf\{t: \sum_{x}H_{2t}(x,x) \le 1+\varepsilon^2 \}=\inf\{t: \sum_{x} \pi(x) d_{2,x}^2(t)  \le \varepsilon^2 \}, \]
where $d_{2,x}^2(t):= \frac{H_{2t}(x,x)}{\pi(x)}-1 \ge 0$ (see \eqref{eq: generalLp}-\eqref{eq: taupeps} below). The $\varepsilon$ $L_{2}$ mixing time, starting from initial state $x$, is  $\mix^{(2),x}(\eps):=\inf\{ t:\frac{H_{2t}(x,x)}{\pi(x)}-1 \le \eps^2 \}$. We write $\mix^{(\infty)}$, $\mix^{(2),x}$ and  $t_{\mathrm{ave-mix}}^{(2)}$ for  $\mix^{(\infty)}(1/2)$,  $\mix^{(2),x}(1/2)$  and $t_{\mathrm{ave-mix}}^{(2)}(1/2)$. We note that \eqref{e:1.2intro1.2} below follows from \eqref{e:1.2intro1.2'}  by considering $x$ such that $t_{\mathrm{mix}}^{(2),x}(\frac{1}{\sqrt{\eps}})=\half t_{\mathrm{mix}}^{(\infty)}(\eps)   $ (such $x$ exists by \eqref{eq: generalLp} below).  
\begin{maintheorem}
\label{p:mixequiv}
For an irreducible reversible Markov chain with a finite state space  we have\footnote{We write $o(1)$ for terms which vanish as $n \to \infty$. We write $f_n=o(g_n)$ or $f_n \ll g_n$ if $f_n/g_n=o(1)$. We write $f_n=O(g_n)$ (and also $g_n=\Omega(f_n)$) if there exists a constant $C>0$ such that $|f_n| \le C |g_n|$ for all $n$. We write  $f_n=\Theta(g_n)$ if  $f_n=O(g_n)$ and  $g_n=O(f_n)$. We use  $\lesssim$ and $\asymp$  to indicate that an inequality (or a two-sided inequality, respectively) holds up to a universal constant for all Markov chains under consideration (usually, all reversible chains with a finite state space, or such chains that are also~transitive).} 
\begin{equation}
\label{e:1.2intro1.2}
\begin{split}
 t_{\mathrm{mix}}^{(\infty)}(\varepsilon) & \le  \rel \max \left\{1,  \log \left(\frac{ \max_{x} t_{\pi \to x}}{\varepsilon \rel} \right) \right\} \text{ for all $\eps >0 $ and so}
\\  t_{\mathrm{mix}}^{(\infty)} & \lesssim \rel \log(1+t_{\mathrm{hit}}/\rel ).
\end{split} 
\end{equation}
\begin{equation}
\label{e:1.2intro1.2'}
\forall \, x \text{ and all }\eps>0, \quad t_{\mathrm{mix}}^{(2),x}(\varepsilon) \le \frac{1}{2} \rel \max \left\{1,  \log \left( \frac{ t_{\pi \to x}}{\eps^2\rel} \right) \right\}. 
\end{equation}
\begin{equation}
\label{e:tavetodot}
 t_{\mathrm{ave-mix}}^{(2)} \le \frac{1}{2} \rel   \log \left(\frac{4 t_{\odot}}{\rel} \right).
 \end{equation}
\end{maintheorem}
It is classical (see \eqref{e:trelmix} and \eqref{e:mixlehit}) that 
\begin{equation}
\label{e:classicthm1}
\rel \log 2 \le \mix^{(\infty)} \le \min \left\{ \rel \log \left( 4/\pi_*  \right),9 \hit \right\}, \qquad \text{where }\pi_*:=\min_x \pi(x).
 \end{equation} In light of the general bound $t_{\mathrm{hit}} \le  \frac{2e}{(e-1)\pi_*}\rel$ (see \eqref{e:stupidboundonthit}) one sees that \eqref{e:1.2intro1.2} refines \eqref{e:classicthm1}. For instance, for simple random walk on a 2 dimensional (discrete grid) torus $\mathbb{Z}_m^2$ of side length $m$ \eqref{e:1.2intro1.2} gives $t_{\mathrm{mix}}^{(\infty)} =O( m^2 \log \log m)$ whereas \eqref{e:classicthm1} gives  $t_{\mathrm{mix}}^{(\infty)} =O( m^2  \log m)$ (we used the fact that in this case $\rel =\Theta( m^2)$ and $\hit = \Theta( m^2 \log m)$, see e.g.\ \cite{LPW}).  

Typically \eqref{e:1.2intro1.2} is much better than the bound $\mix^{(\infty)} \le 9 \hit$ (indeed, $a \log(1+b/a) \le b$ for all $a,b>0$ but  $a \log(1+b/a) \ll b$ when $a \ll b$; apply this to $a=\rel$ and $b=\hit \ge \rel=a$ which typically satisfy $a \ll b$). A  rare\footnote{In a work in preparation with Benjamini, Tessera and Tointon we show that for a sequence of finite vertex-transitive graphs $G_n$ of fixed degree and increasing sizes satisfying that their mixing times are proportional to their  maximal hitting times,  $G_n$ rescaled by their diameters converge in the Gromov-Hausdorff topology to  the unit circle $S^1$. Thus under transitivity $\mix^{(\infty)}  =\Theta( \hit)$ is a very strong condition.} transitive example where this fails is simple random walk on the $n$-cycle in which $\mix^{(\infty)} = \Theta( n^2 )=  \hit = \Theta( \rel \log(1+t_{\mathrm{hit}}/\rel ))$) (without transitivity, it is still fairly  uncommon for  $\mix^{(\infty)}$ and $\hit$ to be of the same order). 

The improvement of  \eqref{e:1.2intro1.2} over $\mix^{(\infty)} \le \rel \log \left( 4 /\pi_*\right)$ is usually less significant than its improvement over  $\mix^{(\infty)} \le 9 \hit$. This is because  $\mix^{(\infty)} \le \rel \log \left( 4 /\pi_*\right)$    is already reasonably sharp, considering the complementary inequality $\mix^{(\infty)} \ge \rel \log 2$. In order for \eqref{e:1.2intro1.2} to improve upon  $\mix^{(\infty)} \le \rel \log \left( 4 /\pi_*  \right)$ by more than a constant factor,  it is necessary that $\log(\hit/\rel ) \ll  \log \left( 1 /\pi_*  \right) $, which for a random walk on an $n$-vertex graph is the same as the condition that $\log( \hit/\rel) \ll \log n$. (This forbids e.g., $\rel \le n^{1-\Omega(1)}$.) 

Other than the 2D torus, there are no many natural examples  satisfying that  $\log(\hit/\rel ) \ll  \log \left( 1 /\pi_*  \right) $ and $\rel \ll \hit$ (especially ones which are transitive). (We mention the condition $\rel \ll \hit$ as the case $\rel = \Theta( \hit)$ was already discussed above).  Theorem \ref{thm:higherorder} provides a higher order versions of  Theorem \ref{p:mixequiv} which can be used to recover the correct order of the mixing time in more examples (see Example \ref{example:higherorder}).

 Equation \eqref{e:41} below gives a probabilistic  interpretation for the r.h.s.\ of \eqref{e:1.2intro1.2} in terms of hitting times for branching random walk. 

Throughout the superscript `$(n)$' indicates that we are considering the $n$th Markov chain in the sequence (not to be confused with `$(\infty)$' and `$(2)$' for $L_{\infty}$ and $L_2$ mixing times). 
\begin{corol}
\label{cor:corofthm1}
For a sequence of irreducible reversible Markov chains with finite state spaces  the following are equivalent:
\begin{itemize}
\item[(i)] $\left(\mix^{(\infty)}\right)^{(n)} =\Theta \left( \hit^{(n)}\right) $ (respectively,  $\left(\mix^{(\infty)}\right)^{(n)} \ll \hit^{(n)} $),
 \item[(ii)]   $\rel^{(n)} =\Theta \left( \hit^{(n)}\right)$ (respectively,  $\rel^{(n)} \ll \hit^{(n)} $).
\end{itemize}
 Also,   $\left(t_{\mathrm{ave-mix}}^{(2)}\right)^{(n)} =\Theta \left( t_{\odot}^{(n)}\right) $ (respectively,  $\left(t_{\mathrm{ave-mix}}^{(2)}\right)^{(n)}\ll t_{\odot}^{(n)} $)
if and only if   $\rel^{(n)} =\Theta \left( t_{\odot}^{(n)}\right)$ (respectively, $ \rel^{(n)} \ll t_{\odot}^{(n)}$). 
\end{corol}
The equivalence between (i) and (ii) is immediate from \eqref{e:1.2intro1.2} in conjunction $\rel  \le   t_{\odot} \le \hit$. Likewise, the second equivalence in the corollary follows at once from \eqref{e:tavetodot} above. 

Let us put Corollary \ref{cor:corofthm1} in a broader context. Spectral conditions play an important role in the modern theory of  Markov chains. A common theme is that for a sequence of reversible Markov chains with finite state spaces of diverging sizes,  certain phenomena can be understood in terms of the simple condition that the product of the spectral-gap and some other natural quantity diverges.  One instance is the cutoff phenomenon and the well-known product condition. This is the condition that the product of the mixing time and the spectral-gap diverges. It is a necessary condition for precutoff in total-variation \cite[Proposition 18.4]{LPW} and  a necessary and sufficient condition for cutoff in $L_2$ \cite{Chen} (in total variation it is known to be a sufficient condition for random walks on trees \cite{basu}; see the recent \cite{Salez} where a modified spectral condition is shown to imply (total variation) cutoff for chains with non-negative curvature). Another such example is given in \cite{hermonspec}, where it is shown that for a sequence of reversible transitive Markov chains (or more generally, ones for which the average and maximal hitting times of states are of the same order) the cover time is concentrated around its mean (for all initial states) if and only if the product of the spectral-gap and the (expected) cover time diverges (this refines a classical result of Aldous \cite{aldouscover}). 

Corollary \ref{cor:corofthm1} concerns the condition that the product of the spectral-gap and the maximal hitting time diverges. This condition was first studied in the context of hitting times in transitive reversible chains by Aldous \cite{aldoushitting}, where it is shown to imply that the maximal and the average hitting time differ only by a smaller order term, and that the law of the hitting time of a vertex is close to an exponential distribution for most initial states.

\textbf{Mean field behavior for coalescing random walks:}  \textcolor{black}{As we now explain, Corollary \ref{cor:corofthm1} resolves a conjecture of Aldous and Fill \cite[Open Problem 14.12]{aldousfill} (re-iterated more recently by Aldous \cite[Open Problem 5]{aldousslide}). The conjecture asserts that for a sequence of vertex-transitive graphs, the condition that $\rel^{(n)} \ll \hit^{(n)} $ implies mean-field behavior for the coalescing time $\tau_{\mathrm{coal}} $ of coalescing random walks. The term mean-field behavior here means that if $t_{\mathrm{meet}}$ is the ``meeting-time", which is defined as the expected collision time of two independent walks started each at equilibrium, then the law of  $\tau_{\mathrm{coal}}/t_{\mathrm{meet}} $ converges in distribution (as the index of the graph diverges) to the corresponding limit for the complete graph on $n$ vertices, which is the law of the coalescence time in Kingman's coalescence (up to a scaling by a factor $n-1$).}

Oliveira \cite{CRW} has already done all of the heavy lifting, and Corollary \ref{cor:corofthm1} merely closes the gap between what is proven in  \cite{CRW} and the aforementioned conjecture. Indeed,  \textcolor{black}{Oliveira \cite{CRW} verified the conjecture (for vertex-transitive graphs) under the seemingly stronger condition  $(\mix^{\TV})^{(n)} \ll \hit^{(n)} $ (see the two comments at the top of p.\  3423 in \cite{CRW}). However, Theorem 3 asserts that this condition is in fact equivalent to the condition  $\rel^{(n)} \ll \hit^{(n)} $. We strongly believe that by combining Oliveira's  methodology from \cite{CRW} with the one from \cite{hermonspec} it is possible to show that for a sequence of reversible chains on finite state spaces $\Omega_n$ (of diverging sizes) with stationary distributions $\pi_n$ satisfying that $\max_{x \in \Omega_n} \pi_{n}(x)=\Theta(\min_{x \in \Omega_n} \pi_{n}(x))$ and $\min_{x \in \Omega_n} \E_{\pi_n}[T_{x}] =\Theta( \max_{x \in \Omega_n} \E_{\pi_n}[T_{x}])  $, the condition that $\rel^{(n)} \ll \hit^{(n)} $ implies mean-field behavior for the coalescing time of coalescing random walks. (Crucially, one can show that such a sequence satisfies $t_{\mathrm{meet}}^{(n)} =\Theta( \hit^{(n)} ) $.)  } 

\textbf{Stability of the condition   $(\mix^{(\infty)})^{(n)} \ll \hit^{(n)}$:} The total-variation mixing time is given by \[\mix^{\TV}:=\inf \left\{t:\max_{x} d_{1,x}(t)\le \frac{1}{2} \right\}, \quad \text{where} \quad d_{1,x}(t):=\sum_{y}|H_{t}(x,y)-\pi(y) | .\]  
While it is known that the order of the $L_{\infty}$ and total variation mixing times can change under a small perturbation of the edge weights or under a rough-isometry \cite{Ding,unifsensitivity,HK,hermonsen}, it follows from Corollary \ref{cor:corofthm1} and \eqref{e:trelmix} that the condition  $(\mix^{(\infty)})^{(n)} \ll \hit^{(n)} $ as well as the condition  $(\mix^{\TV})^{(n)} \ll \hit^{(n)}$ are robust under rough-isometries as well as under a small perturbation of the edge weights, a fact which  a-priori is entirely non-obvious.\footnote{The fact that the maximal (expected) hitting time can change only by a bounded factor under a quasi isometry can be seen from the commute-time identity (e.g.\ \cite[Eq.\ (10.14)]{LPW})   combined with the robustness of the effective-resistance under quasi isometries (cf.\ the proof of Theorem 2.17 in \cite{LP}).} This is in contrast with the spectral condition $\rel^{(n)} \ll (\mix^{\TV})^{(n)}$, which is shown in \cite[Thm.\ 3]{hermonsen} to be sensitive to small perturbations.

\subsection{Higher order analogs of Theorem \ref{p:mixequiv}}
\label{s:higherorder}
We  write $\rho_x:=\rho_{x,2}$ and $\sigma_x=\sigma_{x,2}$ where for $\ell \in \mathbb{N}$  \begin{equation}
\label{e:rhoellsigmaelldef}
\rho_{x,\ell}:= \int_{0}^{2\ell \rel}\frac{s^{\ell-1} \left(H_s(x,x)-\pi(x)\right)}{(\ell-1)! \pi(x)} \mathrm{d}s \; \text{ and} \quad \sigma_{x,\ell}:= \int_{0}^{\infty}\frac{s^{\ell-1} \left(H_s(x,x)-\pi(x)\right)}{(\ell-1)! \pi(x)} \mathrm{d}s.
\end{equation}
We denote the size of the state space by $n$ and write \[\mathcal{Q}=\mathcal{Q}_2 :=\sum_{i=2}^{n} \frac{1}{\lambda_i^2} \quad \text{and} \quad \mathcal{Q}_{\ell} :=\sum_{i=2}^{n} \frac{1}{\lambda_i^{\ell}},\] where as above $0=\lambda_1<\lambda_2 =\mathrm{gap}\le \cdots \le \lambda_n \le 2 $ are the eigenvalues of $I-P$. The motivation behind these definitions shall become clear soon.

\begin{maintheorem}
\label{thm:higherorder}
For an irreducible reversible Markov chain with a finite state space $V$ and $\ell \in \mathbb{N}$ we have that \begin{equation}
\label{e:higherorder}
\forall \, \eps>0, \quad t_{\mathrm{ave-mix}}^{(2)}(\eps) \le \sfrac 12 \rel \max\{ \log \left(\eps^{-2}  \mathcal{Q}_{\ell}/\rel^{\ell} \right), \ell \}. 
\end{equation}
\begin{equation}
\label{e:higherorder2}
\forall \, x \in V, \eps>0, \quad \mix^{(2),x}(\varepsilon) \le \sfrac{1}{2} \rel \max \left\{ \log \left( \varepsilon^{-2}  \sigma_{x,\ell}/\rel^{\ell} \right) ,\ell \right\}     .
\end{equation}
\begin{equation}
\label{e:higherorder3}
\forall \, \eps>0, \quad \mix^{(\infty)}(\eps) \le  \rel \max \left\{ \log \left( \varepsilon^{-1} \sigma_{\ell}  /\rel^{\ell} \right) ,\ell \right\}, \quad \text{where} \quad \sigma_{\ell}:=\max_{x \in V} \sigma_{x,\ell}.   
\end{equation}    
\end{maintheorem}
Recall that the $\Gamma(\ell,1)$ distribution has density function $f_{\ell}(t)=\frac{t^{\ell-1}e^{-t}}{(\ell-1)!} \Ind{t \ge 0}$. For $\ell \in \N$ it is the law of a sum of $\ell$ i.i.d. Exp$(1)$. Denote $\kappa_{\ell}:=\Pr[Z_{\ell} \le 2\ell ] \asymp 1$, where $Z_{\ell} \sim \Gamma(\ell,1) $.
\begin{propo}
\label{propo:higherorderQell}
For an irreducible reversible Markov chain with a finite state space and $\ell \in \mathbb{N}$ we have that
\begin{equation}
\label{e:45}
   \mathcal{Q}_{\ell}=     \sum_x \pi(x) \sigma_{x,\ell}.
 \end{equation}
 \begin{equation}
\label{e:sigmarhocomparison}
\forall \, x, \quad  \kappa_{\ell} \sigma_{x,\ell}\le \rho_{x,\ell} \le  \sigma_{x,\ell}.
 \end{equation}
\end{propo}
\begin{remark}
Using material from \S\ref{s:def} (namely, \eqref{e:alphax}-\eqref{e:RT}) we have that $\mathcal{Q}_1=t_{\odot}$ and that $\sigma_{x,1}=t_{\pi \to x}$. Hence the case $\ell=1$ corresponds to Theorem \ref{p:mixequiv}. The case $\ell=2$ will also play a special role in what comes. 
\end{remark}
\begin{example}
\label{example:higherorder}
The quantity  $\rho_{x,\ell}$ (which in the transitive setup is proportional to $  \mathcal{Q}_{\ell}$ by Proposition \ref{propo:higherorderQell})  can easily be estimated for random walk on a $d$-dimensional torus $\mathbb{Z}_m^d$ of side length $m$. In particular, taking in \eqref{e:higherorder2} $\ell=2d$ gives     $\mix^{(\infty)}  (\mathbb{Z}_m^{d})=O(d\rel )$,  which is sharp up to a constant factor which depends on $d$. Indeed     $\mix^{(\infty)}  (\mathbb{Z}_m^{d})=\Theta(\rel \log d )$. We note that taking $\ell=2$ already suffices to recover the order of the mixing time up to a constant factor for $d=2,3$ and up to a $\log \log m$ factor for $d=4$.

One can reach the same conclusion using \eqref{e:higherorder} via a harder calculation exploiting the fact that the eigenvalues for random walk on  $\mathbb{Z}_m^d$ have a simple explicit expression (e.g.\ \cite[\S\S 12.3-13.4]{LPW}). We omit this calculation as it is somewhat tedious and not particularly illuminating. \end{example}

Let $(X_t)_{t \ge 0}$ and $(Y_t)_{t \ge 0}$ be two independent realizations of the rate 1 continuous-time Markov chain corresponding to transition matrix $P$. The intersection-time is defined as \[\tau_I:=\inf  \{t: X_t \in \{Y_s: s \in [0,t] \} \text{ or }  Y_t \in \{X_s: s \in [0,t] \}  \}.\]  
It is shown in \cite{inter} that under reversibility\footnote{In \cite{inter} discrete-time lazy chains are considered. However their analysis extends to the continuous-time setup.} \[\mix^{\TV} \lesssim t_{\mathrm{I}}:=\max_{x,y}\E_{x,y}[\tau_I ] \quad \text{where} \quad \E_{x,y}[\cdot]:=\E[\, \cdot \mid X_0=x,Y_0=y ], \] and that if in addition $P$ is also transitive then
\begin{equation}
\label{e:tIQ} t_{\mathrm{I}}\asymp \E_{\pi,\pi}[\tau_I ] \asymp \sqrt{ \mathcal{Q}},
\end{equation}
 where $\E_{\pi,\pi} $ is the expectation when $X_0 \sim \pi$ and $Y_0 \sim \pi$ (independently). It is shown in \cite{inter} (Lemma 3.7) that in the transitive reversible setup $\mixin (1/4)\le 2  \sqrt{ \mathcal{Q}}$ (see \S\ref{s:concluding} for a generalization of this fact). The case $\ell=2$ of \eqref{e:higherorder} offers  a substantial improvement in the transitive reversible setup, improving  $2  \sqrt{ \mathcal{Q}}$ to $\mix^{(\infty)}(1/4) \le  \rel \max \left\{ \log \left( 4 \mathcal{Q}  /\rel^{2} \right) ,2 \right\}$.   Theorem \ref{thm:BRWintro} below (namely \eqref{e:46}) offers an analog of \eqref{e:tIQ} by interpreting this upper bound from \eqref{e:higherorder} in terms of the (expected) intersection time for a certain branching random walk, which we describe in the next subsection. 

Using \eqref{e:tIQ},   the following corollary is an immediate consequence of \eqref{e:higherorder}. 

\begin{corol}
\label{cor:BRWcor}
In the above setup, for a sequence of reversible transitive chains:   \[\rel^{(n)} \ll t_{\mathrm{I}}^{(n)} \quad  \text{if and only if} \quad (\mix^{(\infty)} )^{(n)} \ll t_{\mathrm{I}}^{(n)} .\] \end{corol} 
 
\subsection{Hitting and intersection times for ``critical" branching random walk}
\label{s:interintro}
A branching random walk (\textbf{BRW}) with rate $\gamma$ (think of $\gamma$ as the spectral-gap) is a continuous time process in which each particle splits into two particles at rate $\gamma$, independently of the rest of the particles.\footnote{We can treat other variants, including working in discrete-time and/or having a random offspring distribution supported on $\N:=\{1,2,\ldots\}$. The important thing is that the number of particles at time $t$ grows like $\exp(\Theta(t \, \mathrm{gap} )) $. See Remark \ref{rem:mu} for more details.} Each particle performs a rate 1 continuous-time random walk corresponding to some transition matrix $P$, which we assume to be reversible w.r.t.\ $\pi$, independently of the rest of the particles.      
We consider the case that initially there is a single particle whose initial distribution is the stationary distribution $\pi$.

Theorem \ref{p:mixequiv}  explores the relation between hitting-times, mixing-times and the relaxation-time. The hitting-times of a single state are often much larger than the mixing-time (mixing-times are in fact equivalent to hitting-times of large sets \cite{aldous1982some,basu,PS,Olivehit,hermontech,L2}), and so it is interesting to relate  hitting times of a BRW with $\gamma=\mathrm{gap}$ to mixing-times. As we now explain, the choice  $\gamma=\mathrm{gap}$ is natural. With this choice, the number of particles grow by a constant factor every $1/\mathrm{gap} $ time units. The analog of $1-\mathrm{gap} $ for infinite irreducible reversible chains on a countable state space  is the spectral-radius $\rho$ (see e.g., \cite[\S6.2]{LP}). It is classical that $\rho \le 1$ and that when $\rho= 1$ a branching random walk with offspring distribution of mean $\mu>1$ is recurrent, while when $\rho< 1$ the critical mean offspring distribution for recurrence  of a branching random walk (when the process is conditioned to survive) is $\mu_{\mathrm{c}}=1/\rho $ (e.g.\ \cite{criticalBRW}). At  $\mu_{\mathrm{c}} $  the number of particle grows by a constant factor every $1 / (1-\rho) $ time units. Since  $1/\mathrm{gap} $ is the finite setup analog of   $1 / (1-\rho) $, we may interpret our BRW as a ``critical BRW". It is thus less surprising that such a BRW has interesting connections with the mixing time of the chain.    

Let $\mathcal{T}_x $ be the first time at which state $x$ is visited by a particle. The \emph{intersection-time} of two independent BRWs as above, with independent initial distributions, sampled from $ \pi$, denoted by $\mathcal{I}$,    is defined to be the first time $t$ at which a particle from one of the two processes visits a state which was previously visited by a particle from the other process.

The following two theorems and  Corollary \ref{cor:BRWcor}   refine Theorem \ref{p:mixequiv} and Corollary \ref{cor:corofthm1}. We emphasize that the implicit constants below are all independent of the choice of $P$.  We recall that $t_{\pi \to x}=\E_{\pi}[T_x]$ is the expected hitting time of $x$ by the chain (i.e., by a single particle), starting at equilibrium.  

\begin{maintheorem}
\label{thm:BRWintro}
In the above setup, with $\gamma$ taken to be the spectral-gap of $P$, we have that
\begin{equation}
\label{e:41}
\mix^{(\infty)} \lesssim \rel \log(1+ \hit/\rel)\asymp \max_x \E_{\mathrm{BRW}}[\mathcal{T}_x], \end{equation} 
\begin{equation}
\label{e:42}
\begin{split}
\forall \, x, \quad &    \rel \log\left(t_{\pi \to x}/\rel \right)\Ind{t_{\pi \to x} \ge C _{0}\rel }  \lesssim \E_{\mathrm{BRW}}[\mathcal{T}_x]\lesssim \mix^{\TV} + \rel \log \left(1+t_{\pi \to x}/\rel \right), 
\end{split}
\end{equation}
for some absolute constant $C_0$ (independent of $P$ and $x$).
\end{maintheorem}
Let $\rho_{\mathrm{max}}:=\max_x \rho_x$ and $\rho_{\mathrm{min}}:=\min_x \rho_x$, where  $\rho_x=\rho_{x,2}$ is as in \eqref{e:rhoellsigmaelldef}.
\begin{maintheorem}
\label{thm:BRWintro2}
 In the above setup, with $\gamma$ taken to be the spectral-gap of $P$, we have that
\begin{equation}
\label{e:43}
\E_{\mathrm{BRW}}[\mathcal{I}] \lesssim  \rel \log \left(1+ \frac{\sqrt{ \rho_{\mathrm{max}} }}{\rel}  \right),  
 \end{equation}
 \begin{equation}
\label{e:43b}
 \E_{\mathrm{BRW}}[\mathcal{I}] \gtrsim \rel \log \left(1+  \frac{\sqrt{ \rho_{\mathrm{min}} }}{ \rel}\right) \Ind{\rho_{\mathrm{min}} \ge C  _{1}\rel^2}  ,  
 \end{equation}
for some absolute constant $C_1$ (independent of $P$).
Finally, if $P$ is also transitive then
\begin{equation}
\label{e:46}
\mix^{(\infty)} \le \rel \max\{1,\sfrac 12 \log \left(4  \mathcal{Q}/\rel^{2} \right) \} \asymp \rel \log \left(1+ \sqrt{ \mathcal{Q}}/\rel \right)\asymp  \E_{\mathrm{BRW}}[\mathcal{I}] .
\end{equation}
\end{maintheorem}
We note that the second `$\asymp$'   in  \eqref{e:46} follows from \eqref{e:43}-\eqref{e:43b} using \eqref{e:45} and \eqref{e:sigmarhocomparison}, whereas the inequality in  \eqref{e:46}  is immediate from \eqref{e:higherorder}.

We note that \eqref{e:46} refined \eqref{e:41} in the transitive setup. Indeed,  
it is not hard to see that $\E_{\mathrm{BRW}}[\mathcal{I}] \lesssim \max_x \E_{\mathrm{BRW}}[\mathcal{T}_x]$. Moreover, it is also immediate from $\mathcal{Q}=\sum_{i=2}^{n} \frac{1}{\lambda_i^2}$ (by definition) and the fact that $t_{\odot}=\sum_{i=2}^{n} \frac{1}{\lambda_i}$ (see \eqref{e:RT} below) that $\sqrt{ \mathcal{Q}} \le t_{\odot} \le \hit $.  
\begin{remark}
\label{rem:mu}
We note that we could have assumed that at  rate $\gamma \ge \mathrm{gap}$  each particle splits into a random number of particles with  mean $\mu \ge 1 $ such that $\mu-1 \asymp  \mathrm{gap}/\gamma $ and with a finite second moment $\hat \mu$. The above bounds still hold, with the implicit constants depending only on  $(\mu-1)\gamma/  \mathrm{gap}$ and  $\hat \mu$. Similarly, we could have worked in discrete-time and made the offspring distribution $\nu$ of each particle be supported on $\N:=\{1,2,\ldots\} $. In this setup, at each step each particle makes a step according to $P$ (independently of the rest of the particles), then gives birth to a random number of offspring (with law $\nu$, independently of the rest of the particles) and then vanishes. If the mean of $\nu$ is $1+ \Theta( \mathrm{gap}) $ then the same bounds as above hold (up to a constant factor), with $\rel$ replaced by the absolute relaxation-time, which is defined as $\max\{ \frac{1}{1-|\la|}:\la \neq 1 \text{ is an eigenvalue of }P \}$.   
\end{remark}
\begin{remark}
\label{rem:wostinitial}
It is natural to consider the case where for $\mathcal{T}_x $ the starting point of the BRW is a worst-case starting state, rather than stationary. Likewise, for $\mathcal{I}$ it is natural to consider the case that the two BRWs start from a worst pair of initial states. It is easy to reduce the setup of worst-case starting point(s) to the setup of stationary starting point(s), by starting with a burn-in period of duration $\Omega( \mix^{\TV}) $. Indeed, by the following proposition the upper bounds in Theorems \ref{thm:BRWintro} and \ref{thm:BRWintro2} are all $\Omega(\mix^{\mathrm{TV}})$, so allowing such a burn-in period does not increase their order.
\end{remark}
We believe that $  \E_{\mathrm{BRW}}[\mathcal{I}] \asymp \rel \log(1+ \sqrt{ \mathcal{Q}}/\rel) $ whenever $\mathcal{Q} \asymp \rho_{\mathrm{max}}$. A weaker statement that appears to not require much additional work is that this holds whenever
 $\rho_{\mathrm{min}} \asymp \rho_{\mathrm{max}}$.


\subsection{Our spectral optimization methodology}    
Our methodology involves a `spectral optimization' technique, which to the best of our knowledge is novel. For simplicity, let us first describe the method in the transitive setup. In this case, it is well-known (see \S\ref{s:mix} below) that the $L_{\infty}$ mixing time is given by $\inf \{t:\sum_{i=2}^n e^{-\lambda_i t} \le 1/2 \}$ where as above $0=\lambda_1 \le \lambda_2 \le \cdots \le \lambda_n \le 2$ are the eigenvalues of the Laplacian $I-P$.
There are other
 quantities associated with a reversible Markov chain that can be expressed or bounded in terms of the eigenvalues of $I-P$. One such quantity is the average hitting time which satisfies  $t_{\odot}= \sum_{i=2}^n 1/\lambda_i  $ (see \eqref{e:RT}). The idea is to make an informed choice for $t$, which serves as a candidate upper bound on the $L_{\infty}$  mixing time, and consider the following optimization problem, with parameteres $\beta_2,\beta_3,\ldots,\beta_n \in [\lambda_2 ,\infty]$,  \[\max \left\{ \sum_{i=2}^n e^{-\beta_i t}: \beta_2,\beta_3,\ldots,\beta_n \in [ \lambda_2,\infty] \right\} \quad \text{subject to} \quad  \sum_{i=2}^n 1/\beta_i =t_{\odot}. \]
To avoid  rounding issues, we relax the problem and consider
  \[\max \left\{ \sum_{i=2}^n a_{i}e^{-\beta_i t}: \beta_2,\beta_3,\ldots,\beta_n \in [ \lambda_2,\infty],a_{2 },\ldots,a_n \ge 0 \right\} \quad \text{subject to} \quad  \sum_{i=2}^n a_i/\beta_i =t_{\odot}. \]
To bound the mixing time by the informed choice of $t$ (which we soon make) one only needs to establish that the solution to the above optimization problem is at most $1/2$. In practice, ignoring rounding issues, our analysis shows that  whenever  $\lambda_2 \ge \frac{1}{t}$, for the original optimization problem    the worst case is when $\lambda_2 t_{\odot}$ of the $\beta_i$'s equal  $\lambda_2 $ and the rest of them equal $\infty$. Our proofs will not ignore such rounding issues (effectively, by working with the relaxed optimization problem) and will also address the existence of a maximizer, but for the sake of simplicity we ignore these issues in the current discussion. This allows us to make an informed choice of $t$ by considering $t_*=\rel\log \left(\frac{2 t_{\odot}}{\rel} \right)$ which solves $\lambda_2 t_{\odot}e^{-\lambda_2t_{*}}=\frac 12$ and taking  $t=\max\{ \rel,t_*\} \ge \frac{1}{\lambda_2}$.

To see why the maximum is attained when $\beta_2=\cdots=\beta_{\lambda_2 t_{\odot}+1}=\lambda_2$ and $\beta_j=\infty$ for $j>\lambda_2 t_{\odot}+1$ we note that since $h(x)=xe^{-xt}$ is decreasing for  $x \ge 1/t$ we have that $\lambda_2e^{-\lambda_2t} \ge \beta_i e^{-\beta_it} $ for all $i$ for $t \ge 1/\lambda_2$. This means that replacing some $(a_{i},\beta_i)$ for $i$ such that $\beta_i$ is strictly larger than $\la_2$ by $(a_i\lambda_2/\beta_i,\lambda_2)$ preserves the constraint $\sum_{i=2}^n a_{i}/\beta_i =t_{\odot}$ and increases the value of $ \sum_{i=2}^n a_{i}e^{-\beta_i t}$.    

 In fact, using $t_{\mathrm{ave-mix}}^{(2)}=\inf \{t:\sum_{i=2}^n e^{-2\lambda_i t} \le 1/4 \}$,  the above is also a sketch proof of the inequality \eqref{e:tavetodot} (without assuming transitivity; since here we have $e^{-2\lambda_i t}$ instead of $e^{-\lambda_i t} $, we need that $2t \ge \rel$ in order that  $\lambda_2e^{-2\lambda_2t} \ge \beta_i e^{-2\beta_it} $ for all $i$, instead of $t \ge \rel$ as we needed above).   

The above sketch of the method was done in the transitive setup. In the general case the uniform mixing time of the chain starting from a state $x$ can be expressed as  \[\inf \{t:\sum_{i=2}^n c_i(x) e^{-\lambda_i t} \le 1/2 \},\] where $c_2(x),\ldots,c_n(x)$ are certain coefficient coming from the spectral-decomposition (e.g., \cite[\S12.1]{LPW}). The key is that by \eqref{e:alphax} the expected hitting time of $x$ when the chain starts from its stationary distribution $\pi$, denoted by $t_{\pi \to x}$,  equals $\sum_{i=2}^n c_{i}(x)/\lambda_i$ (crucially, with the same coefficients!). The relevant optimization problem becomes
 \[\sup \left\{ \sum_{i=2}^n a_{i}e^{-\beta_i t}: \beta_2,\ldots,\beta_n \in[ \lambda_2,\infty] \text{ and }a_2,\ldots,a_n \ge 0 \right\} \; \text{ subject to} \quad  \sum_{i=2}^n \frac{a_{i}}{\beta_i} = t_{\pi \to x}. \]
Our analysis shows that whenever  $\lambda_2 \ge \frac{1}{t}$   the supremum is attained with $(a_2,\beta_2)=(\lambda_2 t_{\pi \to x},\lambda_2)$ and $(a_i,\beta_i)=(0,\infty)$ for all $i \ge 3$. Here the relevant choice of $t $ is $\max\{\rel,s_{\ast} \}$, where $s_{\ast}$ is the solution to $\lambda_2  t_{\pi \to x}e^{-\lambda_2s_*}=\frac 12$.       

 A similar analysis can be carried out for  other quantities which can be expressed or bounded in terms of the eigenvalues of $I-P$. This is how we prove Theorem \ref{thm:higherorder}.

\textbf{Sharpness of our bounds:} We conclude the discussion of the method by discussing its sharpness. We already gave a few examples (namely, Example  \ref{example:higherorder} and the discussion after Theorem \ref{p:mixequiv}). Our bounds only involve $\rel$ and $t_{\odot}$ or $\hit$. It is natural to ask if one can derive an improved bound on the mixing time which only involves these quantities. This question is intimately related to the question ``how wasteful is it to assume that all eigenvalues of $I-P$ are equal to $\lambda_2$?", which we effectively do when solving the relevant optimization problem. More precisely, the optimal solution above was $\beta_i=\lambda_2$ for all $i$ such that $\beta_i \neq \infty$. In practice, this can correspond to a certain number $k$ (the above solution corresponds to  $k=\lambda_2 t_{\odot}$) of the eigenvalues of $I-P$ which are extremely close to $\lambda_2$ and the rest are ``too large to matter".

Typically the second eigenvalue does not have such a high multiplicity as $\lambda_2 t_{\odot}$, which results in \eqref{e:1.2intro1.2}-\eqref{e:tavetodot} typically not being sharp (other than when $\frac{1}{\lambda_2}$ is of the same order as $t_{\odot}$ or $\hit$). Other than Ramanujan graphs, it is hard to think of cases of $n$-vertex graphs such that random walk on them satisfies that  $\frac{1}{\lambda_2} \ll \hit$ and such that the second smallest eigenvalue of $I-P$ has multiplicity $n^{\Omega(1)}$ (or that there are at least  $n^{\Omega(1)}$ eigenvalues which are  close to $\la_2$).

On the other hand, one can construct such examples which are birth and death chains. The following construction is due to Ding, Lubsetzky and Peres  \cite[bottom of p.\ 77]{DLPbd}. Consider a birth and death chain on $\{1,\ldots,n\}$ with transition probabilities
\begin{itemize}
\item
$P(i,i)=1-\lambda$ for $2 \le i <n$, $P(1,1)=1-\lambda(1-\varepsilon) $ and  $P(n,n)=1-\lambda\varepsilon $,     
\item
$P(i,i-1)=\lambda \eps$ for $2 \le i \le n$, and  
\item
$P(i,i+1)=\lambda (1- \eps)$ for $1 \le i <n$.  
\end{itemize}
We allow $\la$ to depend on $n$ and pick $\varepsilon=\varepsilon(n,\lambda)$ to vanish sufficiently rapidly in terms of $n$ and $\lambda$. It is not hard to check that this can be done so that all $n$ eigenvalues of $P$ are arbitrarily close to $\lambda$ (even when $\lambda$ depends on $n$) in the sense that they equal $\la(1 \pm o(1))$. If we wish for just $k=k(n)$ of the eigenvalues to be arbitrarily close to $\la$ (for some arbitrary $k<n$) and for the rest of the eigenvalues to be arbitrarily close to $1/2$ we can replace $\la$ by $1/2$ in $P(i,j)$ (where $j=i,i\pm 1$) for all $i \le n-k$. We omit the details.

\subsection{Organization of the paper}
In \S\ref{s:def} we introduce notation and definitions and prove Proposition \ref{propo:higherorderQell}. In \S\ref{s:thm3} we prove Theorems \ref{p:mixequiv} and \ref{thm:higherorder}. In \S\ref{s:inter} we prove Theorems \ref{thm:BRWintro} and \ref{thm:BRWintro2}.

\section{Preliminaries and notation}
\label{s:def} 
 Let $(X_t)_{t=0}^{\infty}$ be an irreducible reversible  Markov chain on a finite state space $V$ with transition matrix $P$ and stationary distribution $\pi$. Denote the law of the continuous-time rate 1 version of the chain starting from vertex $x$ (resp.\ initial distribution $\mu$) by $\mathbb{P}_x$ (respectively, $\Pr_{\mu}$). Denote the corresponding expectation by $\mathbb{E}_x$ (respectively, $\E_{\mu}$). For further background on mixing and  hitting times  see \cite{aldousfill,LPW}.

\subsection{Mixing times and $L_p$ norms}
\label{s:mix} 
 The $L_p$ norm and variance of a function $f \in \R^{V}$ are $\|f\|_p:=(\mathbb{E}_{\pi}[|f|^{p}])^{1/p}$ for $1 \le p < \infty$ (where $\mathbb{E}_{\pi}[h]:= \sum_x \pi(x)h(x)$ for $h \in \R^{V}$) and  $\|f\|_{\infty}:=\max_x |f(x)|$. For $p \in [1,\infty]$ the $L_p$ norm of a signed measure $\sigma$ on $ V$ is
\begin{equation*}
\label{eq: Lpdef}
\|\sigma \|_{p,\pi}:=\|\sigma / \pi \|_p, \quad \text{where} \quad (\sigma / \pi)(x):=\sigma(x) / \pi(x).
\end{equation*}
We denote the worst-case $L_p$ distance at time $t$   by $d_{p}(t):=\max_x d_{p,x}(t)$, where $d_{p,x}(t):= \|\Pr_x^t-\pi \|_{p,\pi}$.  Under reversibility for all $x \in V$ and $t \ge 0 $ (e.g.\ \cite[Prop.\ 4.15]{LPW}) we have that
\begin{equation}
\label{eq: generalLp}
d_{2,x}^2(t)= \sfrac{H_{2t}(x,x)}{\pi(x)}-1 \quad \text{and} \quad d_{\infty}(t):=\max_{x,y} |\sfrac{H_{t}(x,y)}{\pi(y)}-1|=\max_{y} \sfrac{H_{t}(y,y)}{\pi(y)}-1.
\end{equation}
The $\epsilon$ $L_{p}$ \emph{mixing time} of the chain (respectively, for initial state $x$) is defined as
\begin{equation}
\label{eq: taupeps}
\begin{split}
\mix^{(p)}(\epsilon)&:= \min \{t:  d_{p}(t)\le \epsilon \}.
\\ t_{\mathrm{mix}}^{(p),x}(\eps)&:= \min \{t:  d_{p,x}(t)\le \epsilon \}.
\end{split}
\end{equation}
When $\epsilon=1/2$  we omit it from the above
notation. The $\eps$ \emph{total variation mixing time} is defined as $\mix^{\TV}(\eps):=\mix^{(1)}(2\eps) $. We write $\mix^{\TV}:=\mix^{(1)}(\half) $.\footnote{Recall that the total-variation distance is $\|\mu-\nu\|_{\mathrm{TV}}:=\half \|\mu-\nu\|_{1,\pi} $.}   Clearly, $\mix^{(p)} $ is non-decreasing in $p$. Finally, we define the \emph{average $\eps $ $L_2$ mixing time} as 
\begin{equation}
\label{e:aveL2}
\begin{split}
& t_{\mathrm{ave-mix}}^{(2)}(\eps):=\min \left\{t: \sum_{v \in V }\pi(v) d_{2,v}^2(t) \le \eps^2 \right\}
\\ & \text{(by \eqref{eq: generalLp})} =\min \left\{t: \sum_{v \in V }H_{2t}(v,v)\le 1+\eps^2 \right\}=\min \left\{t: \sum_{i=2}^{|V|}e^{-2\la_i t} \le \eps^2 \right\}.
\end{split}
\end{equation}
We now recall  
 the hierarchy between the various quantities considered above. Under reversibility we have that (e.g.\ \cite[Theorems 10.22, 12.4, 12.5 and Lemmas 4.18 and 20.5]{LPW}) 
\begin{equation}
\label{e:trelmix}
\forall \, \eps \in (0,1), \quad \sfrac{1}{\la_2} |\log \eps | \le \mix^{\TV} ( \eps/2) \le  \mix^{(2)} (\eps)=\half  \mixin (\eps^2) \le \sfrac{1}{\la_2} |\log( \eps^{2} \min_x \pi(x)  )|,
\end{equation}
\begin{equation}
\label{e:mixlehit}
 \mixin \le 9 \hit. 
\end{equation}

\subsection{Hitting-times}
\label{s:hitbackground}
We now present some background on hitting times. The random target identity (e.g.~\cite[Lemma 10.1]{LPW}) asserts that $t_{\odot}:= \sum_y \pi(y) \mathbb{E}_x[T_y] $ is independent of $x$, and that under reversibility, for all $x \in V$, writing $a_i(x):=P^{2i}(x,x)+P^{2i+1}(x,x)-2\pi(x)$, we have that   
\begin{equation}
\label{e:alphax}
 t_{\pi \to x} :=  \E_{\pi}[T_x]=\frac{1}{\pi(x)} \sum_{i=0}^{\infty}a_{i}(x)=\frac{1}{\pi(x)} \int_{0}^{\infty}\left( H_t(x,x)-\pi(x) \right)\mathrm{d}t
\end{equation}

(see e.g.~\cite[Proposition 10.26]{LPW} or \cite[Ch.\ 2]{aldousfill}). In fact, this holds even without reversibility, if we define $a_i(x):=-p\pi(x)+ \sum_{j=0}^{p-1}P^{pi+j}(x,x)$, where $p$ is the period of the chain. Averaging over $x$ yields (in the reversible setup) the eigentime identity (see e.g.\ \cite[Proposition 3.13]{aldousfill})  
\begin{equation}
\label{e:RT}
t_{\odot}=\sum_{x,y} \pi(x)\pi(y) \mathbb{E}_x[T_y]=\sum_{y} \sum_{i=0}^{\infty}(P^{i}(y,y)-\pi(y))=\sum_{i=0}^{\infty}\mathrm{[Trace}(P^i)-1]=\sum_{i \ge 2}\frac{1}{\la_i}. \end{equation}
Let $U \sim \pi $ be independent of the chain. Noting that $T_x \le T_U+\inf\{t:X_{t+T_U}=x \}$ and using the random target lemma to argue $\E[T_{U}]=t_{\odot} $, as well as the strong Markov property to argue that $\E[\inf\{t:X_{t+T_U}=x \}]=\E_{\pi}[T_x]=t_{\pi \to x} $, yields: 
\begin{fact}[\cite{LPW} Lemma 10.2]
\label{f:a}
  $\max_x t_{\pi \to x}  \le \hit \le t_{\odot}+\max_x t_{\pi \to x} \le  2\max_x t_{\pi \to x} $.
\end{fact}
The following material can be found at \cite[\S 3.5]{aldousfill}. Some of the relevant ideas are due to the pioneer work of Keilson on total positivity \cite{Keilson}. Under reversibility, for any set $A$ the law of its hitting time $T_A:=\inf\{t:X_t \in A\} $ under initial distribution $\pi$ conditioned on $A^{\complement} $, is a mixture of Exponential distributions, whose minimal parameter $\lambda(A)$ is the Dirichlet eigenvalue of the set $A^{\complement}$. There exists a distribution $\mu_A$, known as the \emph{quasi-stationary distribution} of $A^{\complement}$, under which  $T_A$ has an Exponential distribution of parameter $\la(A)$.  It follows that $\la(A) \ge \frac{1}{\max_a \E_a[T_A]} $. We see that for all $t \ge 0$,  \begin{equation}
\label{e:exptailpi} \Pr_{\pi}[T_y> t] \le \exp(-t/\hit) , \quad \text{and so} \quad \E_{\pi}[T_y^2] \le 2 \hit^{2}. 
\end{equation} 
Using the above description of the law of $T_A$ it is not hard to show (e.g.\ \cite[p.\ 86]{aldousfill}) that
\begin{equation}
\label{e:poscorfortail} \forall \, s,t \ge 0, \qquad \Pr_{\pi}[T_y> t+s \mid T_y \ge s ] \ge \Pr_{\pi}[T_y> t ]. 
\end{equation} 
It follows from the spectral decomposition (e.g., \cite[\S12.1]{LPW}) that for all $x$ and  all $s,t \ge 0$ we have that   
 \begin{equation}
 \label{e:specd}
0< H_{t+s}(x,x) -\pi(x) \le e^{-s/\rel}(H_{t}(x,x) -\pi(x)). 
\end{equation} This easily implies the following lemma: 

\begin{lemma}
\label{lem:auxexpdecay} For every irreducible, reversible Markov chain on a finite state space with a stationary distribution $\pi$, for every state $x$ and all $M > 0$ we have that
\begin{equation}
\label{e:auxexp1}
\int_0^{\infty}(H_{s}(x,x)-\pi(x))\mathrm{d}s \le \frac{e^{M}}{e^M-1}    \int_0^{M\rel}(H_{s}(x,x)-\pi(x))\mathrm{d}s.
\end{equation}
\begin{equation}
\label{e:auxexp2}
\int_0^{\infty}s(H_{s}(x,x)-\pi(x))\mathrm{d}s \le \left(1+\sum_{i=1}^{\infty}(i+1)e^{-iM} \right)   \int_0^{M\rel}s(H_{s}(x,x)-\pi(x))\mathrm{d}s.
\end{equation}
\end{lemma}
\emph{Proof.}
By \eqref{e:specd}
 $\frac{ f(i)}{ f(0) } \le e^{-Mi} $, where $f(i):=\int_{iM \rel }^{(i+1)M\rel}(H_s(x,x)-\pi(x)) \mathrm{d}s$. This easily implies \eqref{e:auxexp1}. Likewise writing   $g(i):=\int_{iM \rel }^{(i+1)M\rel}s(H_s(x,x)-\pi(x)) \mathrm{d}s$,  \eqref{e:auxexp2} follows from $\frac{ g(i)}{ g(1) } \le (i+1) e^{-Mi} $ (which again follows from \eqref{e:specd}).  \qed  

\vspace{3mm}

Combined with Fact \ref{f:a}, \eqref{e:alphax} and \eqref{e:auxexp1}  yield the following weak bound which we record solely for the sake of the discussion in the introduction. Recall that $\pi_*=\min_x \pi(x)$.
\begin{corol}
\label{cor:stupidboundonthit}
For every irreducible, reversible Markov chain on a finite state space $\Omega$ with a stationary distribution $\pi$, we have that
\begin{equation}
\label{e:stupidboundonthit}
\hit \le 2\max_x t_{\pi \to x} \le \frac{2e}{e-1}\frac{1}{\pi_{*}}\int_0^{\rel}(H_{s}(x,x)-\pi(x))\mathrm{d}s \le \frac{2e}{e-1}\rel\frac{1-\pi_{*}}{\pi_{*}} . 
\end{equation}
\end{corol}
\noindent \textbf{Proof of Proposition \eqref{propo:higherorderQell}.} 
We first prove \eqref{e:45}. Let $\ell \in \mathbb{N}$. Using $\int_0^{\infty}\frac{x^{\ell} s^{\ell-1}}{(\ell-1)!}e^{-sx}\mathrm{d}s=1$ for all $x>0$ and the spectral decomposition (as well as a change of order of summation and integration) we have that
\[\mathcal{Q}_{\ell} = \sum_{i=2}^{n}  \int_0^{\infty}\frac{s^{\ell-1}}{(\ell-1)!}e^{-s \lambda_i}\mathrm{d}s=  \int_0^{\infty}\frac{s^{\ell-1}}{(\ell-1)!}\left(\mathrm{Trace}(H_{s})-1\right)\mathrm{d}s=\sum_x  \pi(x) \sigma_{x,\ell}. \] 
We now prove \eqref{e:sigmarhocomparison}. The inequality $ \rho_{x,\ell} \le  \sigma_{x,\ell}$ is trivial. We now prove that  $\kappa_{\ell} \sigma_{x,\ell} \le  \rho_{x,\ell}$. By the spectral decomposition $\frac{H_s(x,x)-\pi(x)}{\pi(x)}=\sum_{i=2}^n f_i(x)^2e^{-\lambda_is}$ where $f_1=1,f_{2},\ldots,f_n$ are the orthonormal basis  (w.r.t.\ the inner-product $\langle f,g\rangle_{\pi}:=\sum_{y}\pi(y)f(y)g(y)$)  corresponding to the eigenvalues $0=\lambda_1 \le \lambda_2 \le \cdots$  $\le \lambda_n \le 2 $ of $I-P$ (i.e.\ $(I-P)f_i=\lambda_i f_i$ for all $i$). Then
\[\rho_{x,\ell}= \int_{0}^{2\ell \rel}\frac{s^{\ell-1} \left(H_s(x,x)-\pi(x)\right)}{(\ell-1)! \pi(x)} \mathrm{d}s=\sum_{i=2}^n \frac{f_i(x)^2}{\lambda_i^{\ell}} \int_{0}^{2\ell \rel}\frac{\lambda_i^{\ell}s^{\ell-1} }{(\ell-1)!} e^{-\lambda_is} \mathrm{d}s \ge \kappa_{\ell}\sum_{i=2}^n \frac{f_i(x)^2}{\lambda_i^{\ell}},  \] where, recalling that $\kappa_{\ell}=\mathbb{P}[\Gamma(\ell,1) \le 2 \ell]$, the last inequality follows  from  \[\forall \,i \ge 2, \quad \int_{0}^{2\ell \rel}\frac{\lambda_i^{\ell}s^{\ell-1} }{(\ell-1)!} e^{-\lambda_is} \mathrm{d}s=^{} \int_{0}^{2\ell \rel \lambda_i }\frac{s^{\ell-1} }{(\ell-1)!} e^{-s} \mathrm{d}s \ge \kappa_{\ell}, \]
where in the inequality we used the fact   $\lambda_i \ge1/ \rel$. 
Finally
\begin{equation}
\label{sigmaxellviaev}
\begin{split}
\sigma_{x,\ell}& = \int_{0}^{\infty}\frac{s^{\ell-1} \left(H_s(x,x)-\pi(x)\right)}{(\ell-1)! \pi(x)} \mathrm{d}s=\sum_{i=2}^n \frac{f_i(x)^2}{\lambda_i^{\ell}} \int_{0}^{\infty}\frac{s^{\ell-1} }{(\ell-1)!}\lambda_i ^{\ell}e^{-\lambda_is} \mathrm{d}s \\ & =\sum_{i=2}^n \frac{f_i(x)^2}{\lambda_i^{\ell}}. \qquad \qed
\end{split}
\end{equation}

\section{Proof of Theorems \ref{p:mixequiv} and \ref{thm:higherorder}}
\label{s:thm3}
We begin the section with a proof of Theorem \ref{thm:higherorder}. 

\textbf{Proof of Theorem \ref{thm:higherorder}:}
We first prove \eqref{e:higherorder}. Let $t:=\frac 12 \rel \max\{  \log \left(4  \mathcal{Q_{\ell}}/\rel^{\ell} \right),\ell \}$.  Recall that  $ \sum_{v \in V }\pi(v) d_{2,v}^2(t)=\sum_{i=2}^n e^{-2 \la_i t} $ and that  $\mathcal{Q}_{\ell}=\sum_{i=2}^{n} \lambda_i^{-\ell}$. Hence  $ \sum_{v \in V }\pi(v) d_{2,v}^2(t)$ is  bounded by the value of the solution to the optimization problem:
\[\max \sum_{i=1}^{n-1} a_{i}e^{-2\beta_i t}, \]
subject to the conditions 
\begin{itemize}
\item[(1)]  $\sum_{i=1}^{n-1} a_{i}\beta_i^{-\ell} =\mathcal{Q}_{\ell} $,
\item[(2)] $0 \le a_i \le \mathcal{Q}_{\ell}\lambda_n^{\ell} $ for all $i$, and
\item[(3)] $  \la_2 \le \beta_1 \le \beta_2 \le \cdots \le \beta_{n-1}  \le \infty  $. Moreover $\beta_i \le \lambda_n$  for all $i$ such that $\beta_i<\infty$.
\end{itemize}
If we further require that exactly $k$ of the $\beta_i$'s are finite, then a maximum is indeed attained by compactness. As there are at most $n-1$ possible values for $k$ we see that also without such a constraint a maximum must be attained. We will show  that the maximum is attained  when \[\beta_1=\la_2, \; \, a_{1} = \mathcal{Q}_{\ell} \la_2^{\ell}, \quad \text{and} \quad \beta_i=\infty \text{ for all }i \ge 2.  \] 
Thus $ \sum_{v \in V }\pi(v) d_{2,v}^2(t) \le \mathcal{Q}_{\ell} \la_2^{\ell}\exp(-2\la_2 t ) $. Substituting the value of $t$ we see that this is at most $1/4$ as desired. It remains only to verify the above claim that the maximum of the optimization problem is attained when $\beta_1=\lambda_2$ and $a_{1} =\mathcal{Q}_{\ell} \la_2^{\ell}$, while $\beta_i=\infty$ for all $i \ge 2$.  

While it is possible to conclude the proof even without this observation, we note that by a simple Lagrange multipliers calculation one gets that for any maximizer if $\beta_i,\beta_j \notin \{\la_2,\lambda_n,\infty \}$ and $\min \{a_i,a_j\}>0$ then $\beta_i=\beta_j$ (as can be seen by considering the derivatives w.r.t.\ $a_i$ and $a_j$). 

 By collecting together $\beta_{i},\ldots,\beta_{i+j}$ which all equal one another, and increasing the value of $a_i$ to $\sum_{k=0}^{j}a_{i+k}$   we see that we can add to constraint (3) the requirement that $\beta_i<\beta_j$ whenever $i<j$ and $\beta_j<\infty$. Together with the previous paragraph, this  means that the maximum is attained by a solution with one of the following forms:
\begin{itemize}
\item[Case 1] $\beta_1=\lambda_2$,  $a_1=\mathcal{Q}_{\ell}\lambda_2^{\ell}$  and $\beta_i = \infty$ for all $i \ge 2$. 
\item[Case 2] $\beta_1=\lambda_2$, $\lambda_2<\beta_2  \le \lambda_n $, $a_1,a_2>0$ and $\beta_i = \infty$ for all $i \ge 3$. 
\item[Case 3] $\beta_1=\lambda_2$, $\lambda_2<\beta_2 < \lambda_n$, $\beta_3=\lambda_n$,  $a_1,a_2,a_3>0$  and $\beta_i = \infty$ for all $i \ge 4$. 
\item[Case 4] $\lambda_2<\beta_1< \lambda_n $, $\beta_2=\lambda_n$,  $a_1,a_2>0$  and $\beta_i = \infty$ for all $i \ge 3$.
\item[Case 5] $\lambda_2<\beta_1 \le \lambda_n$  $a_1=\mathcal{Q}_{\ell}\beta_1^{\ell}$  and $\beta_i = \infty$ for all $i \ge 2$.  
\end{itemize}
We now rule out cases 2 to 5.         
 In cases 2 to 5 above there is some $i$ such that $a_i=c>0$ and $\beta_i=\mu$ for some  $\lambda_2<\mu \le \lambda_n $ and $c>0$. Replacing $a_i=c$ with $a_{i}=c \lambda_2^{\ell}/\mu^{\ell}$ and $\beta_i=\mu$ with $\beta_i=\lambda_2$ preserves  constraints (1)-(3) (possibly, up to relabeling of the indices). We will show that this solution has a strictly larger value. This means that the maximum cannot be attained by a solution having one of the forms of cases 2 to 5. 

 To check that the above transformation gives rise to a solution with a larger value it suffices to check that   $ce^{-2 \mu t} \le c\frac{ \lambda_2^{\ell}}{\mu^{\ell}} e^{-2 \lambda_2 t}$ which is the same as    $\mu^{\ell}e^{-2 \mu t} \le \lambda_2^{\ell} e^{-2 \lambda_2 t}$. Indeed, by taking derivative w.r.t.\ $x$ we see that the function $h(x)=x^{\ell} e^{-2x t}$ is non-increasing if $\frac{\ell}{2x} \le t $. In other words, since by definition $\frac{\ell}{2\la_2} \le t $ we indeed have that    $\mu^{\ell} e^{-2 \mu t} \le \lambda_2^{\ell} e^{-2 \lambda_2 t}$ for all $\mu \ge \lambda_2$, as desired. This concludes the proof  of \eqref{e:higherorder}.

\noindent \textbf{Proof of \eqref{e:higherorder2}:}  Let $t:= \sfrac{1}{2} \rel \max \left\{ \log \left( \varepsilon^{-2}  \sigma_{x,\ell}/\rel^{\ell} \right) ,\ell \right\}     $.
By the spectral decomposition and \eqref{eq: generalLp}  $d_{2,x}^2(t)=\frac{H_{2t}(x,x)-\pi(x)}{\pi(x)}=\sum_{i=2}^n f_{i}(x)^{2}e^{-2 \la_i t}$, where $f_i$ are an orthonormal (w.r.t.\ the inner-product $\langle f,g\rangle_{\pi}:=\sum_{y}\pi(y)f(y)g(y)$) basis of eigenfunctions with $(I-P)f_i=\lambda_if_i$,  and that by \eqref{sigmaxellviaev}  $\sigma_{x,\ell}=\sum_{i=2}^n \frac{f_i(x)^2}{\lambda_i^{\ell}}$. Hence, similarly as in the proof of \eqref{e:higherorder},   $d_{2,x}^2(t)$ is  bounded by the value of the solution to the optimization problem:
\[\max \sum_{i=1}^{n-1} a_{i}e^{-2\beta_i t}, \]
subject to the conditions
 
\begin{itemize}
\item[(1)]  $\sum_{i=1}^{n-1} a_{i}\beta_i^{-\ell} =\sigma_{x,\ell} $,
\item[(2)] $0 \le a_i \le \sigma_{x,\ell}\lambda_n^{\ell} $ for all $i$, and
\item[(3)] $  \la_2 \le \beta_1 \le \beta_2 \le \cdots \le \beta_{n-1}  \le \infty  $. Moreover $\beta_i \le \lambda_n$  for all $i$ such that $\beta_i<\infty$.
\end{itemize}
The existence of the maximum is justified exactly as in the proof of \eqref{e:higherorder}. Similarly, as in \eqref{e:higherorder} we have that the maximum is attained at a solution with one of the following forms:

\begin{itemize}
\item[Case 1] $\beta_1=\lambda_2$,  $a_1=\sigma_{x,\ell}\lambda_2^{\ell}$  and $\beta_i = \infty$ for all $i \ge 2$. 
\item[Case 2] $\beta_1=\lambda_2$, $\lambda_2<\beta_2  \le \lambda_n $, $a_1,a_2>0$ and $\beta_i = \infty$ for all $i \ge 3$. 
\item[Case 3] $\beta_1=\lambda_2$, $\lambda_2<\beta_2 < \lambda_n$, $\beta_3=\lambda_n$,  $a_1,a_2,a_3>0$  and $\beta_i = \infty$ for all $i \ge 4$. 
\item[Case 4] $\lambda_2<\beta_1< \lambda_n $, $\beta_2=\lambda_n$,  $a_1,a_2>0$  and $\beta_i = \infty$ for all $i \ge 3$.
\item[Case 5] $\lambda_2<\beta_1 \le \lambda_n$  $a_1=\sigma_{x,\ell}\beta_1^{\ell}$  and $\beta_i = \infty$ for all $i \ge 2$.  
\end{itemize}
We now rule out cases 2 to 5. This is done in exactly the same fashion as in the proof of \eqref{e:higherorder}. We include the details for the sake of completeness. 

 In cases 2 to 5 above there is some $i$ such that $a_i=c>0$ and $\beta_i=\mu$ for some  $\lambda_2<\mu \le \lambda_n $ and $c>0$. Replacing $a_i=c$ with $a_{i}=c \lambda_2^{\ell}/\mu^{\ell}$ and $\beta_i=\mu$ with $\beta_i=\lambda_2$ preserves  constraints (1)-(3) (possibly, up to relabeling of the indices). We will show that this solution has a strictly larger value. This means that the maximum cannot be attained by a solution having one of the forms of cases 2 to 5. 
 
To check that the above transformation gives rise to a solution with a larger value it suffices to check that   $ce^{-2 \mu t} \le c\frac{ \lambda_2^{\ell}}{\mu^{\ell}} e^{-2 \lambda_2 t}$ which is the same as    $\mu^{\ell}e^{-2 \mu t} \le \lambda_2^{\ell} e^{-2 \lambda_2 t}$. We already saw above that this is satisfied whenever $\frac{\ell}{2 \lambda_2} \le t $ (which indeed holds by the definition of $t$). The proof  of \eqref{e:higherorder2} is now concluded by noting that by the choice of $t$ (using Case 1) we have that $d_{2,x}^2(t) \le \sigma_{x,\ell}\lambda_2^{\ell}  e^{-2 \lambda_2 t} \le 1/4$.

\noindent \textbf{Proof of \eqref{e:higherorder3}:} this is immediate from \eqref{e:higherorder2} together with the equality $\mix^{(2)} (\eps)=\half  \mixin (\eps^2)$ from \eqref{e:trelmix} (which follows directly from \eqref{eq: generalLp}).      \qed

\vspace{5pt}

\noindent{\textbf{Proof of Theorem \ref{p:mixequiv}:}} This is just the $\ell=1$ case of Theorem \ref{thm:higherorder} by noting that $t_{\odot}=\sum_{i=2}^n 1/\lambda_i=\mathcal{Q}_1$ by \eqref{e:RT} and that $t_{\pi \to x}=\sigma_{x,1}$ by \eqref{e:alphax}. \qed

\section{Hitting and intersection times for branching random walk - Proof of Theorems \ref{thm:BRWintro} and \ref{thm:BRWintro2}}
\label{s:inter}
\emph{Proof of Theorem \ref{thm:BRWintro}:}
We first note that the first inequality in \eqref{e:41}  is a repetition of  \eqref{e:1.2intro1.2}. We now argue that the  relation $\rel \log(1+ \hit/\rel)\asymp \max_x \E_{\mathrm{BRW}}[\mathcal{T}_x]$ in \eqref{e:41} follows from \eqref{e:42}:   $J_{x}\Ind{t_{\pi \to x} \ge C_0  \rel  } \lesssim  \E_{\mathrm{BRW}}[\mathcal{T}_x]\lesssim \mix^{\TV} + J_{x}$, where throughout \[J_x:=\rel \log(1+t_{\pi \to x} /\rel)\] and (as always) $t_{\pi \to x}:=\E_{\pi}[T_x]$.  Indeed, if there exists $x$ such that $t_{\pi \to x} \ge C_0 \rel $, then (if $C_0$ is sufficiently large) by Theorem \ref{p:mixequiv} we can also find  $x$ such that \[ \rel \log(t_{\pi \to x} /\rel)  \ge   \mix^{\TV}/32.   \] For such $x$ by \eqref{e:42} we have that $\E_{\mathrm{BRW}}[\mathcal{T}_x] \asymp J_x \gtrsim \mix^{\TV}  $, from which it is easy to see that $\max_x \E_{\mathrm{BRW}}[\mathcal{T}_x] \asymp \max_{x} J_x \asymp \rel \log(1+\hit/\rel) $ (as $\max_{x}t_{\pi \to x} \asymp \hit $ by Fact \ref{f:a}). 

If   $\max_x t_{\pi \to x} < C_0 \rel $, then again using  $\max_{x}t_{\pi \to x} \asymp \hit $ we see that $\hit \asymp \rel $, and so $\rel \log(1+ \hit/\rel) \asymp \hit $, and so $\max_x \E_{\mathrm{BRW}}[\mathcal{T}_x] \le \max_{x}t_{\pi \to x}   \le  \hit \asymp \rel \log(1+ \hit/\rel)   $ (the first inequality follows by considering the hitting time of $x$ by the first particle). Moreover, if $\hit \asymp \rel $ then with a positive probability we have that the BRW has a single particle by time $t_{\pi \to x}/2 $, where $x$ is picked so that $t_{\pi \to x}=\max_{z}t_{\pi \to z}$. By the Paley-Zygmund inequality we have that $\Pr_{\pi}[T_x> t_{\pi \to x}/2] \ge \frac{1}{4} \frac{t_{\pi \to x}^2}{\E_{\pi}[T_x^2]} \ge \frac{1}{32}$, where we have used $t_{\pi \to x} =\max_{z}t_{\pi \to z}\ge \half \hit $ (by the choice of $x$ and Fact \ref{f:a}) and $\E_{\pi}[T_x^2] \le 2 \hit^2 $ \eqref{e:exptailpi}. 

It follows that when $\max_z   t_{\pi \to z} < C_0 \rel $ we have that for some $x$ such that  $t_{\pi \to x}=\max_z   t_{\pi \to z}  $  with probability  bounded from below (uniformly in $P$) $\mathcal{T}_x \ge \frac{ t_{\pi \to x}}{2} \gtrsim \hit  \asymp \rel \log(1+ \hit/\rel)$. This concludes the derivation of \eqref{e:41} from \eqref{e:42}.   

\textbf{Proof of upper bound in} \eqref{e:42}: We now prove that $  \E_{\mathrm{BRW}}[\mathcal{T}_x]\lesssim \mix^{\TV} + J_{x}$. If $t_{\pi \to x} < e \rel $ then we are done, since in this case  $\E_{\mathrm{BRW}}[\mathcal{T}_x] \le t_{\pi \to x} \asymp J_{x} $. Now assume that  $t_{\pi \to x} \ge e \rel $.

\begin{lemma}
\label{lem:positiveprobab}
Consider the BRW from above, where initially there is a single particle starting from initial position $z$. Denote the corresponding law by $\P_{\mathrm{BRW},z}$. When the initial distribution of the particle is $\pi$ we simply write  $\P_{\mathrm{BRW}}$ for the law and   $\E_{\mathrm{BRW}}$ for the corresponding expectation.  
Assume that for some $M \ge 1$ and $p \in (0,1)$  we have that
\begin{equation}
\label{e:geomtailbrw}
\max_z \P_{\mathrm{BRW},z}[\mathcal{T}_x> 10  \mix^{\TV} +MJ_{x})] \le 1- p. 
\end{equation}
Then 
\begin{equation}
\label{e:geomtailbrw2}
 \E_{\mathrm{BRW}}[\mathcal{T}_x] \le \max_z \E_{\mathrm{BRW},z}[\mathcal{T}_x] \le p^{-1}(10  \mix^{\TV} +MJ_{x}). 
\end{equation}
Moreover, if for some $M \ge 1$,  $K\in \mathbb{N}$ and $p_1 \in (0,1)$ we have that 
\[\P_{\pi}[T_x > MJ_x/2] \le (1- p_1)^{1/K} \quad \text{and} \] \[ \P_{\mathrm{BRW}}[\text{the number of particles at time } MJ_x/2 \text{ is at least }2K] \ge p_2, \]
then \eqref{e:geomtailbrw} holds for this $M$ with $p:= p_1p_2 p_3 $, where $p_3:=\mathbb{P}[\mathrm{Bin}(2K,3/4)  \ge K]$.  
\end{lemma}
\textbf{Proof of Lemma \ref{lem:positiveprobab}:} It is easy to see that \eqref{e:geomtailbrw2} is implied by the inequality
\[\forall \, i,\quad  \max_z \P_{\mathrm{BRW},z}[\mathcal{T}_x> (i+1)(10  \mix^{\TV} +MJ_{x}) \mid \mathcal{T}_x> i(10  \mix^{\TV} +MJ_{x})  ] \le 1- p \]
(indeed, this implies that $\lceil \mathcal{T}_x/(10  \mix^{\TV} +MJ_{x})\rceil $ is stochastically dominated by the Geometric distribution with parameter $p$).
This indeed follows from \eqref{e:geomtailbrw} using the Markov property, by noting that if at time $i(10  \mix^{\TV} +MJ_{x})$ we remove all but one of the particles this can only increase $\mathcal{T}_x$. It remains to prove that  \eqref{e:geomtailbrw} holds with $p:= p_{1}p_2 p_3 $. 

As one might expect from the statement of the lemma, the idea is to first let the process run for $MJ_x/2 $ times units, during which with probability at least $p_2$ the number of particles will become at least $2K$. Let us work on this event (this accounts for the term $p_2$ in $p_1p_2p_3$). We then let these particles run for a burn-in period of length $10 \mix^{\TV}$.
We will argue rigorously that at a cost of the term $p_3$ we may assume that at least $K$ of these particles have i.i.d.\ stationary positions at time $MJ_x/2+10 \mix^{\TV} $ (this accounts for the term $p_3$ in $p_1p_2p_3$).
  Finally, once we have $K$ such particles with i.i.d.\ stationary positions, by letting them run for additional $MJ_x/2$ time units, we see that the probability that at least one of them hits $x$ in these additional $MJ_x/2$ time units is at least
\[1-(\P_{\pi}[T_x > MJ_x/2])^{K} \ge p_1. \]

By the above discussion, to conclude the proof of the lemma, it remains only to show that for any initial state $z$ we can write the law of the walk at time $t=10  \mix^{\TV}$ as a mixture of the form $\frac{3}{4}\pi+\frac 14 \mu_z$ for some distribution $\mu_z$, or equivalently, that $H_t(z,y) \ge \frac 34\pi(y)$ for all $y$. Indeed, such a mixture can be sampled by flipping a bias coin with heads probability $3/4$. If it lands on heads, sample $\pi$ and otherwise sample $\mu_z$. Given $2K$ particles, we can do so independently for each of them (where the coin flips as well as the samples of $\pi$ and $\mu_z$ are all independent, and where for each particle, $z$ corresponds to its position at the beginning of the burn-in period). The probability that at least $K$ of the coins will land on heads is exactly $p_3$. (This part of the argument is taken from the proof of Proposition 4.3 in \cite{CRWrev}, where such an argument is carried out in detail.)

To prove this, we will need to recall the notion of separation mixing time and  its connection with the total variation mixing time.
The $\delta$-separation mixing-time is defined as \[t_{\mathrm{sep}}(\delta):=\inf\{t: \min_{x,y} H_{t}(x,y)/\pi(y) \ge 1-\delta \}.\]
For reversible Markov chains we have that $\mixtv(\eps) \le t_{\mathrm{sep}}(\eps) \le 2 \mixtv(\eps/4)  $ for all $\eps \in (0,\sfrac 14) $ (e.g., \cite[Lemmas 6.16 and 6.17]{LPW}).
In particular, using the `sub-multiplicativity' property  $\mix^{\TV}(2^{-i}) \le i  \mix^{\TV} $ (recalling that $\mix^{\TV}$ is a shorthand notation for $\mix^{\TV}(1/4)$; e.g., \cite[(4.34)]{LPW})
\begin{equation}
\label{e:TVsepsubmult}
10  \mix^{\TV} \ge2 \mix^{\TV}(2^{-4}) \ge t_{\mathrm{sep}}(\sfrac 14) \ge\mix^{\TV}    .
\end{equation}
We are done, since   
 $\min_{x,y} H_{t}(x,y)/\pi(y)$ is non-decreasing (e.g.\ \cite[Exercise 6.4]{LPW}). \qed 

At time $8J_x $ the number of particles is at least $2\lceil t_{\pi \to x}/2\rel \rceil $ with probability bounded from below (uniformly in $x$ and $P$).
Hence we may take above $M=16$ and $K=\lceil t_{\pi \to x}/2\rel \rceil$, and take $p_2$ and $p_3$ to be uniformly bounded away from zero (over all $x$ and $P$).
 
 It remains to bound $p_1$ from below. Let $K=\lceil t_{\pi \to x}/2\rel \rceil$ and $M=16$ be as above.  By \eqref{e:poscorfortail} (used in the second inequality) and Markov's inequality we have that \[\left(\P_{\pi}[T_x > MJ_x/2] \right)^K \le (\Pr_{\pi}[T_x > e \rel])^{K} \le \Pr_{\pi}[T_x > e \rel K]  \le \E_{\pi}[T_x ]/(e \rel K) \le 2/e, \]
where in the first inequality we have used the fact that $MJ_x/2=8J_x \ge e \rel$ by our assumption that   $t_{\pi \to x} \ge e \rel $. This concludes the proof of the upper bound in \eqref{e:42}.   

\textbf{Proof of lower bound in} \eqref{e:42}: We now show that $  \E_{\mathrm{BRW}}[\mathcal{T}_x]\gtrsim  J_{x}$ when    $t_{\pi \to x} \ge C_{0}\rel$, provided $C_0$ is sufficiently large.  Write $N_x(t):=\int_0^t \text{(number of particles at $x$ at time }s) \mathrm{d}s $.  Then 
\begin{equation}
\label{NJx}
\P_{\mathrm{BRW}}[\mathcal{T}_x \le \sfrac 14J_{x}]=\P_{\mathrm{BRW}}[N_x(\sfrac 14J_{x})>0] \le \frac{\E_{\mathrm{BRW}}[N_x(\sfrac 12 J_{x})]}{\E_{\mathrm{BRW}}[N_x(\sfrac 12 J_{x}) \mid  N_x(\sfrac 14J_{x})>0]}.
 \end{equation}  
By the Markov property, and using \eqref{e:auxexp1} and the fact that $J_x \ge16 \rel $, provided that $C_{0}$ is sufficiently large, we have that \[\E_{\mathrm{BRW}}[N_x(\sfrac 12 J_{x}) \mid  N_x(\sfrac 14J_{x})>0] \ge \int_{0}^{\sfrac 14J_{x}}H_s(x,x)\mathrm{d}s \ge \frac{7}{8} \int_{0}^{\infty }(H_s(x,x)-\pi(x)) \mathrm{d}s.  \]
Recall that $t_{\pi \to x}=\sfrac{1}{\pi(x)}  \int_{0}^{\infty }(H_s(x,x)-\pi(x)) \mathrm{d}s$ (\eqref{e:alphax}) and so  \[  \E_{\mathrm{BRW}}[N_x(\sfrac 12 J_{x}) \mid  N_x(\sfrac 14J_{x})>0] \ge \sfrac 78 \pi(x)t_{\pi \to x} . \] On the other hand, by stationarity, and using the fact that the expected number of particles at time $t$ is $2^{t/\rel} $, we have that \[\E_{\mathrm{BRW}}[N_x(\sfrac 12 J_{x})]=\pi(x) \int_0^{\sfrac 12 J_{x}}2^{s/\rel}\mathrm{d}s \le \sfrac{ \pi(x)\rel}{ \log 2}2^{J_x/(2 \rel)} \le \pi(x) t_{\pi \to x}/2 ,  \]
where the last inequality follows by the definition of $J_x$, provided $C_{0}$ is sufficiently large. Plugging the last two estimates into \eqref{NJx}, we see that  $\P_{\mathrm{BRW}}[\mathcal{T}_x \le \sfrac 14J_{x}] \le 4/7$ , which implies that  $  \E_{\mathrm{BRW}}[\mathcal{T}_x] \gtrsim  J_{x}$, as desired. \qed

\vspace{5pt}

\emph{Proof of Theorem \ref{thm:BRWintro2}:} 
We first prove \eqref{e:43b}. In fact, we prove a slightly stronger statement. Before stating it we need to introduce a new quantity (we expect that typically $\chi(2) \asymp \rho_{\mathrm{min}} $, recalling that $\rho_{\mathrm{min}}:=\min_x \rho_x$ where  $\rho_x=\rho_{x,2}= \int_{0}^{4 \rel}\frac{s \left(H_s(x,x)-\pi(x)\right)}{ \pi(x)} \mathrm{d}s$ is as in \eqref{e:rhoellsigmaelldef})\footnote{By  the spectral-decomposition $H_s(x,x)$ is decreasing in $s$ (and $H_s(x,x)\searrow \pi(x) $ as $s \to \infty$), which implies that $ \int_{0}^{ t}\frac{sH_s(x,x)}{\pi(x)} \mathrm{d}s \ge \frac{H_{t}(x,x)}{2 \pi(x)}t^2 $. This, together with $H_1(x,x) \ge 1/e $  implies that for some absolute constant $c_0>0$, for all $a>1$ we have that $ \chi_{x}(a/2) \ge \max \{ 2 \mix^{(2),x}(\sqrt{a-1}),c_0 (a \pi(x))^{-1/2} \}$.}:  \[\chi(a):=\min_x \chi_{x}(a), \quad \text{where} \quad \chi_{x}(a):=\inf \left\{t>0: \int_{0}^{ t}\frac{sH_s(x,x)}{\pi(x)} \mathrm{d}s \le a t^2 \right\}. \]
 \[L:= \begin{cases} c\rel \log \left(1+\frac{\sqrt{ \rho_{\mathrm{min}} }}{\rel} \right)   & \text{if } \rho_{\mathrm{min}} \ge e^{8/c}  \rel^{2} \text{ or } \chi(2^{10}) \ge \rel  \\
\chi(2^{10}) & \text{if } \rho_{\mathrm{min}} < e^{8/c}    \rel^{2} \text{ and } \chi(2^{10})<\rel \\
\end{cases} ,\] for some absolute constant $c \in (0,1)$ to be determined later. The lower bound in  \eqref{e:43b} is an immediate consequence of the following proposition together with the definition of $L$.  

\begin{propo}
\label{p:lowerBRW}
For some absolute constant $c \in (0,1)$ (from the definition of $L$) we have that   $L \lesssim \E_{\mathrm{BRW}}[\mathcal{I}]$. 
\end{propo} 

\textbf{Proof of Proposition \ref{p:lowerBRW}:} We first observe that it suffice to show that $\P_{\mathrm{BRW}}[\mathcal{I} \le L ]$ is bounded away from 1.

 We write $U(x,s)$ (respectively, $V(x,s) $) for  the number of particles from the first (respectively, second) BRW which occupy state $x$ at time $s$. Consider \[M(t):=\sum_x \int_0^t\int_0^t \frac{U(x,s)V(x,r)}{\pi(x)} \mathrm{d}s \mathrm{d}r .\]  Then 
\begin{equation}
\label{M2L}
\P_{\mathrm{BRW}}[\mathcal{I} \le L ]=\P_{\mathrm{BRW}}[M( L)>0] \le \frac{\E_{\mathrm{BRW}}[M(2L)]}{\E_{\mathrm{BRW}}[M(2L) \mid  M(L)>0]}.
\end{equation}  
By the Markov property and the independence between the two BRWs,  we have that
\begin{equation}
\label{M2L2}
\E_{\mathrm{BRW}}[M(2L) \mid  M(L)>0] \ge \min_x \int_0^{L}\int_0^{L}\sum_y \frac{H_s(x,y)H_r(x,y)}{\pi(y)}\mathrm{d}s\mathrm{d}r.
 \end{equation}
To see this, observe that if $M(L)>0$ there is some $x$ which is occupied at some time $s _{0}\le L$ by a particle from the first BRW and at some time $r_{0} \le L$ by a particle from the second BRW. Considering the contribution to $M(2L)$ of these two particles, corresponding to times $s \in [s_0,2L] \supset [s_0,s_0+L]$ for the first particle and to times $r\in [r_0,2L] \supset [r_0,r_0+L]$ for the second particle yields the r.h.s.\ of \eqref{M2L2}. 
 
By reversibility $\sum_y \frac{H_s(x,y)H_r(x,y)}{\pi(y)}=H_{s+r}(x,x) $ and so \eqref{M2L2} 
 \begin{equation}
 \label{e:53ML}
  \E_{\mathrm{BRW}}[M(2L) \mid  M(L)>0] \ge \min_{x}  \int_{0}^{L}\frac{sH_s(x,x)}{\pi(x)}\mathrm{d}s.
 \end{equation}
On the other hand, by stationarity and independence of the two BRWs, and using the fact that for each of the BRWs the expected number of particles at time $t$ is $2^{t/\rel} $, we have that
 \begin{equation}
 \label{e:53ML2}
\E_{\mathrm{BRW}}[M(2L)]=\sum_x \pi(x) \int_0^{L}\int_0^{L}2^{(s+r)/\rel}\mathrm{d}s \mathrm{d}r= \sfrac{ \rel^{2}}{( \log 2)^{2}}\left( 2^{L/ \rel} -1 \right)^2.
 \end{equation} 
\textbf{Case I -  $\rho_{\mathrm{min}} \ge e^{2/c}  \rel^{2}$: }First consider the case that $\rho_{\mathrm{min}} \ge e^{8/c}  \rel^{2}$. Then $L \ge 4\rel $, and so by the definition of $\rho_{\mathrm{min}} $ we have that  \[\E_{\mathrm{BRW}}[M(2L) \mid  M(L)>0] \ge\min_{x}  \int_{0}^{L}\frac{sH_s(x,x)}{\pi(x)}\mathrm{d}s \ge \rho_{\mathrm{min}}  .\]
By \eqref{e:53ML2}, if $c$ is sufficiently small we have that
\[\E_{\mathrm{BRW}}[M(2L)] \le   \rho_{\mathrm{min}}/2. \]
Plugging the last two estimates into \eqref{M2L} concludes the proof in this case.

 \textbf{Case II -  
$\rho_{\mathrm{min}} < e^{8/c}   \rel^{2}$ and $\chi(2^{10}) < \rel $:} Then $L =  \chi(2^{10})  $.  Using the fact that $H_s(x,x)$ is decreasing in $s$, and using the definition of $\chi(\cdot)$, we see that
\begin{equation}
\label{e:48210}
\min_{x}  \int_{0}^{L}\frac{sH_s(x,x)}{\pi(x)}\mathrm{d}s \ge \min_{x}  \int_{0}^{\chi(2^{10})}\frac{sH_s(x,x)}{\pi(x)}\mathrm{d}s= 2^{10}\chi(2^{10})^2=2^{10}L^{2}  .
\end{equation}

By \eqref{e:53ML2} and the fact that $L =  \chi(2^{10}) < \rel   $ (together with $2^x-1 \le2 x$ doe $0 \le x \le 1$ and $\log 2 > 1/2$),  we have that
\[\E_{\mathrm{BRW}}[M(2L)] \le 16  L^2. \]  
Plugging this, together with the estimate $\E_{\mathrm{BRW}}[M(2L) \mid  M(L)>0] \ge2^{10}L^{2}$ (which follows by combining \eqref{e:48210} and \eqref{e:53ML}) into    \eqref{M2L} concludes the proof in this case.

 \textbf{Case III -  
$\rho_{\mathrm{min}} < e^{8/c}   \rel^{2}$ and $\chi(2^{10}) \ge \rel $:} Finally, consider the case that
$\rho_{\mathrm{min}} < e^{8/c} \rel^{2}$ and $\chi(2^{10}) \ge \rel $.  Then $4L \le \rel \le \chi(2^{10})  $. Hence \[\min_{x}  \int_{0}^{\rel}\frac{sH_s(x,x)}{\pi(x)}\mathrm{d}s \ge \min_{x}  \int_{0}^{4L}\frac{sH_s(x,x)}{\pi(x)}\mathrm{d}s \ge 2^{10}(4L)^2=2^{14}L^2.\]  Using $4L \le \rel $ again, by \eqref{e:53ML2} we have that $\E_{\mathrm{BRW}}[M(2L)] \le \sfrac{ \rel^{2}}{( \log 2)^{2}}\left( 2^{L/ \rel} -1 \right)^2 \le \sfrac{ 1}{(2 \log 2)^{2}}  L^2 $. The proof is hence concluded using \eqref{M2L}. \qed

\vspace{5pt}

\textbf{Proof of }$\E_{\mathrm{BRW}}[\mathcal{I}] \lesssim  \rel \log \left(1+\frac{\sqrt{ \rho_{\mathrm{max}} }}{\rel} \right)$\textbf{:} Recall that $\rho_{\mathrm{max}}=\max_x \rho_x$, where $\rho_{x}=\rho_{x,2}$, and that $\mathcal{Q}=\mathcal{Q}_2$. By Proposition \ref{propo:higherorderQell} $\rel^2 \le \mathcal{Q} \asymp \sum_x \pi(x) \rho_x$, and so  $\sqrt{ \rho_{\mathrm{max}} } \gtrsim \rel $. Hence at time $t=C_{2}\rel \log \left(1+\frac{\sqrt{ \rho_{\mathrm{max}} }}{\rel} \right)$ both BRWs have at least $2C_3( \mathrm{gap}\sqrt{ \rho_{\mathrm{max}} })^2 $ particles with  probability bounded from below (uniformly in $P$), provided $C_2$ is sufficiently large in terms of $C_3$ (where $C_3$ is to be determined). As in the proof  of Lemma \ref{lem:positiveprobab} (which is part of the proof of $  \E_{\mathrm{BRW}}[\mathcal{T}_x]\lesssim \mix^{\TV} + J_{x}$), by taking a burn-in period of duration $t_\mathrm{sep}(1/4)$, we may assume that each of BRWs has $S:=\lceil C_3 ( \mathrm{gap}\sqrt{ \rho_{\mathrm{max}} })^2 \rceil$ particles whose positions are at equilibrium, independently (see 2.\ below). We note for later usage that by Theorem \ref{thm:higherorder} and \eqref{e:TVsepsubmult} $ t_\mathrm{sep}(1/4) \lesssim  \mix^{\TV} \lesssim \rel \log \left(1+\frac{\sqrt{ \rho_{\mathrm{max}} }}{\rel} \right)$ . Hence the duration of  such a burn-in period can be absorbed into the implicit constant.

We can now label the particles in the two processes by $1,\ldots,S:=\lceil C_3 ( \mathrm{gap}\sqrt{ \rho_{\mathrm{max}} })^2 \rceil $ and $1', \ldots ,S'$, respectively. Denote the time interval $[t+t_{\mathrm{sep}}(1/4),t+t_{\mathrm{sep}}(1/4)+\rel]$ by $\mathcal{J}$. As we explain below in more detail, using  the independence of the pairs $(i,i')$,  by only considering intersections of particle $i$ with particle $i'$, it suffices to show that,  
\begin{equation}
\label{e:Ilrel}
\P[\text{the pair of particles }(i,i') \text{ intersect during }\mathcal{J}] =\P_{\pi,\pi}[\tau_I< \rel ] \gtrsim \frac{1}{(\mathrm{gap})^2\rho_{\mathrm{max}} },  
\end{equation}    
where as in  \S\ref{s:interintro}, $\tau_I$ is the intersection time of two independent realizations $(X_t)_{t \ge 0 }$ and $(Y_t)_{t \ge 0 }$ of the Markov chain (with no branching), and $\P_{\pi,\pi}$ indicates that $X_0 \sim \pi$ and  $Y_0 \sim \pi$, independently. Indeed, one way for $\mathcal{I}$ to be at most $t+t_\mathrm{sep}(1/4)+\rel  $ is the following:
\begin{itemize}
\item[1.] At time $t=C_{2}\rel \log \left(1+\mathrm{gap}\sqrt{ \rho_{\mathrm{max}} } \right)$ both BRWs have at least $2C_3( \mathrm{gap}\sqrt{ \rho_{\mathrm{max}} })^2 $ particles. For each $C_3$ we can pick $C_2$ so that this has probability bounded from~below.
\item[2.] For each of the BRWs, after a burn-in period between time $t$ and $t+t_{\mathrm{sep}}(1/4)$ we get at least $S:=\lceil C_3 ( \mathrm{gap}\sqrt{ \rho_{\mathrm{max}} })^2 \rceil$ particles with i.i.d.\ stationary positions. As in the proof Lemma \ref{lem:positiveprobab}, by enlarging the probability space if necessary this can be made precise, and the probability of this is bounded from below.
\item[3.] With the above labeling of the particles, for some $1 \le i \le S$ we have that the paths of particles $i$ and $i'$ intersect during the time interval $\mathcal{J}=[t+t_{\mathrm{sep}}(1/4),t+t_{\mathrm{sep}}(1/4)+\rel]$. By this we mean that for some $r,s \in\mathcal{J}$ and some state $z$ we have that particle $i$ is at $z$ at time $r$   and  particle $i'$ is at $z$ at time $s$. Note that for each such pair, the probability this happens is $\P_{\pi,\pi}[I< \rel ]$ and so the probability this happens for at least one pair $(i,i')$ out of the $S$ pairs is
\[1-\P_{\pi,\pi}[\tau_I \ge \rel ]^S \] 
\end{itemize}
By adjusting $C_2$ from the definition of $t$ we may ensure that the constant $C_3$ from the definition of $S$ can be chosen to be arbitrarily large, while maintaining the probabilities from 1.\ and 2.\ above being bounded away from zero. 

If \eqref{e:Ilrel} holds, then the probability of the event in 3.\ above is bounded away from zero, since  it follows from \eqref{e:Ilrel} that for some constant $c>0$ we have \[\P_{\pi,\pi}[\tau_I \ge \rel ]^S \le \left( 1-  \frac{c}{\mathrm{gap}^2\rho_{\mathrm{max}} }\right)^S \le \exp \left(-\frac{cS}{\mathrm{gap}^2\rho_{\mathrm{max}} } \right)=\exp \left(-\frac{c\lceil C_3 ( \mathrm{gap}\sqrt{ \rho_{\mathrm{max}} })^2 \rceil}{\mathrm{gap}^2\rho_{\mathrm{max}} } \right) \] which can indeed be made arbitrarily small by picking $C_3$ to be large enough in terms of~$\frac 1c$. 

 This  only implies that with probability bounded away from 0 we have that \[\mathcal{I} \le t+t_{\mathrm{sep}}(1/4)+\rel \lesssim  \rel \log \left(1+\mathrm{gap}\sqrt{ \rho_{\mathrm{max}} } \right).\]
(We have already explained above that $t_{\mathrm{sep}}(1/4) \lesssim t \asymp \rel \log \left(1+\mathrm{gap}\sqrt{ \rho_{\mathrm{max}} } \right)$. Since by \eqref{e:TVsepsubmult}  $\rel \lesssim \mix^{\TV} \le t_{\mathrm{sep}}(1/4)$, this shows that indeed $ t+t_{\mathrm{sep}}(1/4)+\rel \lesssim  \rel \log \left(1+\mathrm{gap}\sqrt{ \rho_{\mathrm{max}} } \right)$.)

 However, applying the Markov property as in the proof of Lemma \ref{lem:positiveprobab}, can upgrade such an estimate into the desired estimate  $\E_{\mathrm{BRW}}[\mathcal{I}] \lesssim  \rel \log \left(1+\mathrm{gap}\sqrt{ \rho_{\mathrm{max}} } \right)$.

It remains to prove \eqref{e:Ilrel}, which we now do. For the sake of clarity of presentation, we phrase this as a lemma.
\begin{lemma}
\label{lem:4.9I}
$\P_{\pi,\pi}[\tau_I< \rel ] \gtrsim \frac{1}{\mathrm{gap}^2\rho_{\mathrm{max}} }$.
\end{lemma}

\textbf{Proof:} Consider $R(t):=\sum_x \frac{1}{\pi(x)}\int_0^t\int_0^t \Ind{X_s=x=Y_r} \mathrm{d}s \mathrm{d}r$. Then $\E_{\pi,\pi}[R(t)]=t^2 $ and
\[R(t)^2 =\sum_{x,y} \frac{1}{\pi(x) \pi(y) }\int_0^t\int_0^t\int_0^t\int_0^t \Ind{X_s=x=Y_r,X_{a}=y=Y_{b}} \mathrm{d}s \mathrm{d}r \mathrm{d}a \mathrm{d}b. \] 
Taking expectation, and using reversibility, we see that 
\[\E_{\pi,\pi}[R(t)^2] \lesssim \sum_{x} \pi(x)t^{2} \int_0^t\int_0^t \sum_y \frac{ H_{a}(x,y)H_b(x,y)}{\pi(y)} \mathrm{d}a \mathrm{d}b \] \[= \sum_{x} \pi(x)t^{2} \int_0^t\int_0^t  \frac{ H_{a+b}(x,x)}{\pi(x)} \mathrm{d}a \mathrm{d}b \le  t^2 \max_{x} \int_0^{2t}  \frac{sH_s(x,x)}{\pi(x)}\mathrm{d}s \lesssim t^2 \max_{x} \int_0^t  \frac{sH_s(x,x)}{\pi(x)}\mathrm{d}s, \]
(where we have used the fact that $H_s(x,x) $ is non-decreasing in $s$ to argue that $ \int_t^{2t}  \frac{sH_s(x,x)}{\pi(x)}\mathrm{d}s$ $ \le 16  \max_{x} \int_{t/2}^{t}  \frac{sH_s(x,x)}{\pi(x)}\mathrm{d}s   $). 
Now, for $t=\rel$ we have that \[\max_{x} \int_0^t  \frac{sH_s(x,x)}{\pi(x)} \le \rho_{\mathrm{max}} + \int_{0}^t s \mathrm{d}s=\rho_{\mathrm{max}}+\half \rel^{2} \lesssim \rho_{\mathrm{max}}   ,\] where we have used the fact that by \eqref{e:45} $\rel \le \sqrt{\mathcal{Q}} \lesssim \sum_x \pi(x) \rho_x \le \rho_{\mathrm{max}} $.
Thus by the Payley-Zygmund inequality
\[ \P_{\pi,\pi}[\tau_I< \rel ]=\P_{\pi,\pi}[R( \rel)>0 ]  \ge \frac{(\E_{\pi,\pi}[R(\rel)])^{2}}{\E_{\pi,\pi}[R(\rel)^2] } \gtrsim \rel^2/\rho_{\mathrm{max}}. \quad \qed \] 

\vspace{5pt}

We now prove that if $P$ is  transitive we have that
\[\mix^{(\infty)} \lesssim \rel \log(1+ \sqrt{ \mathcal{Q}}/\rel)\asymp  \E_{\mathrm{BRW}}[\mathcal{I}] .\]
The first inequality follows from \eqref{e:higherorder}. The inequality  $ \E_{\mathrm{BRW}}[\mathcal{I}] \lesssim  \rel \log(1+ \sqrt{ \mathcal{Q}}/\rel) $ follows from $ \E_{\mathrm{BRW}}[\mathcal{I}] \lesssim  \rel \log(1+ \sqrt{ \rho_{\mathrm{\max}}}/\rel) $ and \eqref{e:45}. 

The proof of  \eqref{e:45} gives that $\sfrac 23 \le \mathcal{Q} / \sum_x \pi(x) \rho_x\le \varkappa  $, where $\varkappa:= 1+\sum_{i=1}^{\infty}(i+1)e^{-i/2}$.  When $\mathcal{Q} \ge C_1 \varkappa \rel^2 $ by  \eqref{e:43b} we have that $\E_{\mathrm{BRW}}[\mathcal{I}] \gtrsim  \rel \log(1+ \sqrt{ \mathcal{Q}}/\rel) $. 

It remains to consider the case that $\mathcal{Q} < C_1 \varkappa \rel^2 $ and to show that also in this case we have that   $\E_{\mathrm{BRW}}[\mathcal{I}] \gtrsim  \rel \log(1+ \sqrt{ \mathcal{Q}}/\rel) $. When    $\mathcal{Q} < C_1 \varkappa \rel^2 $, using   $\mathcal{Q} \ge \rel^2 $ we have that $\rel \log(1+ \sqrt{ \mathcal{Q}}/\rel)\asymp \rel $, and so in order to prove  $\E_{\mathrm{BRW}}[\mathcal{I}] \gtrsim  \rel \log(1+ \sqrt{ \mathcal{Q}}/\rel) $, it suffices to show that in this case for some $c_{1},c_{2} \in (0,1)$ we have that
\begin{equation}
\label{e:c1c2Irel}
\P_{\pi,\pi}[\tau_I > c_{1} \rel ]\ge c_{2}.
 \end{equation}
To see this, observe that with probability bounded away from zero we have that both BRWs have only a single particle by time $c_{1} \rel$. Given this event, the conditional probability that $\mathcal{I}>c_{1} \rel$ is precisely $\P_{\pi,\pi}[\tau_I > c_{1} \rel ]$. Finally, in the considered case $\rel \asymp \rel \log(1+ \sqrt{ \mathcal{Q}}/\rel) $, and so it suffices to show that $\P_{\mathrm{BRW}}[\mathcal{I} > c_{1} \rel ] \gtrsim 1$ in the considered case. 

It remains to prove \eqref{e:c1c2Irel} (assuming  $\mathcal{Q} < C_1 \varkappa \rel^2 $), which we now do.
By \eqref{e:tIQ} (taken from \cite[Lemma 3.7]{inter}) we have that $\mixin (1/4)\le 2  \sqrt{ \mathcal{Q}} \asymp \E_{\pi,\pi}[\tau_I]$. Applying Markov's inequality and the Markov property inductively, yields that $\Pr_{\pi,\pi}\left[\tau_I> j \left( 2 \E_{\pi,\pi}[\tau_I] +\mixin (1/4) \right) \right] $ decays exponentially in $j$ (we omit the details, as they are routine). Hence if  $\mathcal{Q} < C_1 \varkappa \rel^2 $ we have that $\E_{\pi,\pi}[\tau_I^2] \lesssim( \E_{\pi,\pi}[\tau_I])^{2} $, and so by the Paley-Zygmund inequality \[\P_{\pi,\pi}[\tau_I > c_{1} \rel ] \ge \P_{\pi,\pi}[\tau_I >\sfrac 12  \E_{\pi,\pi}[\tau_I]] \ge \frac{( \E_{\pi,\pi}[\tau_I])^{2} }{4\E_{\pi,\pi}[\tau_I^2]}  \ge c_{2}\] where we used the fact that  $\rel \le \sqrt{ \mathcal{Q}} \asymp \E_{\pi,\pi}[\tau_I] $ to find an absolute constant $c_1>0$ such that $c_1 \rel \le \sfrac 12\E_{\pi,\pi}[\tau_I]$ and used  $\E_{\pi,\pi}[\tau_I^2] \lesssim( \E_{\pi,\pi}[\tau_I])^{2} $ to find an absolute constant $c_2>0$ as above. This concludes the proof.   \qed   
\section{Concluding remarks}
\label{s:concluding}
\textbf{Two additional inequalities:} Using the fact $H_t(x,x)$ is decreasing in $t$, together with the definition of $\mix^{(2),x}$ we have that \[\sigma_{x,\ell} \ge \int_{0}^{2\mix^{(2),x}}\frac{s^{\ell-1} \left(H_s(x,x)-\pi(x)\right)}{(\ell-1)! \pi(x)} \mathrm{d}s \ge\int_{0}^{2\mix^{(2),x}}\frac{s^{\ell-1} }{(\ell-1)! 4} \mathrm{d}s=\frac{(2\mix^{(2),x})^{\ell}}{\ell!4},   \]       
\[\text{and so} \quad 2\mix^{(2),x}\le (\ell! 4 \sigma_{x,\ell}  )^{1/\ell} \le 4\ell \sigma_{x,\ell}  ^{1/\ell},   \]
i.e.\ $\mix^{(2),x} \le 2\ell \sigma_{x,\ell}  ^{1/\ell}$. Similarly, using \eqref{e:45}
\[ \mathcal{Q}_{\ell}=\sum_x \pi(x) \sigma_{x,\ell} \ge\sum_x \pi(x) \int_{0}^{2t}\frac{s^{\ell-1} \left(H_{2t}(x,x)-\pi(x)\right)}{(\ell-1)! \pi(x)} \mathrm{d}s =\frac{(2t)^{\ell}}{\ell!}\sum_x \pi(x)d_{2,x}^2(t).   \]  
Taking $t=t_{\mathrm{ave-mix}}^{(2)}$ yields that   $t_{\mathrm{ave-mix}}^{(2)}\le 2\ell \mathcal{Q}_{\ell}^{1/\ell}$ (under transitivity,   $t_{\mathrm{mix}}^{(2)}\le 2 \ell \mathcal{Q}_{\ell}^{1/\ell}$).
\begin{problem}
It would be interesting to find a probabilistic interpretation for the bounds from Theorem \ref{thm:higherorder} for general $\ell$ as we have done (in the following subsection) for $\ell=1,2$ in terms of hitting and intersection times of branching random walks. Likewise, it would be interesting to find a probabilistic interpretation for $\mathcal{Q}_{\ell}$ for $\ell \ge 3$.
\end{problem}

{\bf Acknowledgements:} We thank Tom Hutchcroft and Perla Sousi for helpful discussions. We also thank the referees for useful comments that helped improve the presentation of the paper.


\begin{thebibliography}{BKC}


\bibitem{aldoushitting}
Aldous, D., Hitting times for random walks on vertex-transitive graphs. \emph{Math. Proc. Cambridge Philos. Soc.}  \textbf{106} (1989), no. 1, 179--191. \href{http://www.ams.org/mathscinet-getitem?mr=MR0994089}{\textcolor{blue}{MR0994089}}.


\bibitem{aldousslide}
Aldous, D., Mixing times and hitting times. (2010). Available at \href{http://www.stat.berkeley.
edu/~aldous/Talks/slides.html}{\textcolor{blue}{http://www.stat.berkeley.
edu/~aldous/Talks/slides.html}}


\bibitem{aldous1982some}
Aldous, D.,
\newblock Some inequalities for reversible Markov chains.
\newblock {\em Journal of the London Mathematical Society}, 2(3):564--576,
  1982. \href{http://www.ams.org/mathscinet-getitem?mr=MR657512}{\textcolor{blue}{MR657512}}

\bibitem{aldouscover}
Aldous, D., Threshold limits for cover times. \emph{J. Theoret. Probab.} \textbf{4} (1991), no. 1, 197--211. \href{http://www.ams.org/mathscinet-getitem?mr=MR1088401}{\textcolor{blue}{MR1088401}}.





\bibitem{aldousfill}
Aldous, D., and Fill, J.,
  \emph{Reversible Markov chains and random walks on graphs}. Unfinished manuscript. Available at the first author's \href{https://www.stat.berkeley.edu/~aldous/}{\textcolor{blue}{website}}.








\bibitem{basu}
Basu, R., Hermon, J., and Peres, Y., Characterization of cutoff for reversible Markov chains. \emph{Ann. Probab.} 45 (2017), no. 3, 1448--1487. \href{http://www.ams.org/mathscinet-getitem?mr=MR3650406}{\textcolor{blue}{MR3650406}}



 












\bibitem{Chen}
Chen, G.\ Y., and Saloff-Coste, L., The cutoff phenomenon for ergodic Markov processes.
 \newblock \emph{Electron. J. Probab.} 13 (2008), no. 3, 26--78.  
\href{http://www.ams.org/mathscinet-getitem?mr=MR2375599}{\textcolor{blue}{MR2375599}}


\bibitem{Ding}
Ding, J., and Peres, Y.
\newblock
 Sensitivity of mixing times.
\newblock
{\em Electron. Commun. Probab.} \textbf{18} (2013), paper 88, 6pp. 
\newblock Available at: \href{https://projecteuclid.org/euclid.ecp/1465315627}{\nolinkurl{projecteuclid/1465315627}}
\bibitem{DLPbd}
Ding, J., Lubetzky, E., and Peres, Y.
\newblock
 Total variation cutoff in birth-and-death chains.
\newblock
{\em Probab. Theory Related Fields} \textbf{146} (2010),  no. 1-2, 61--85. 
\newblock   
\href{http://www.ams.org/mathscinet-getitem?mr=MR2550359}{\textcolor{blue}{MR2550359}}


 


 

 


\bibitem{criticalBRW}
Gantert, N., and M\"{u}ller, S.,
\newblock The critical branching {M}arkov chain is transient.
\newblock \emph{Markov Process. Related Fields} 12 (2006), no. 4, 805--814. \href{http://www.ams.org/mathscinet-getitem?mr=MR2284404}{\textcolor{blue}{MR2284404}}

\bibitem{hermonspec}
Hermon, J.,
\newblock A spectral characterization for concentration of the
cover time.
\newblock (2019) To appear in \emph{Journal of Theoretical Probab.} Arxiv preprint \href{https://arxiv.org/abs/1809.00145}{\textcolor{blue}{arXiv:1809.00145}}


\bibitem{hermontech}
Hermon, J.,
\newblock A technical report on hitting times, mixing and cutoff.
\newblock \emph{ALEA Lat. Am. J. Probab. Math. Stat.} 15 (2018), no. 1, 101--120. \href{http://www.ams.org/mathscinet-getitem?mr=MR3765366}{\textcolor{blue}{MR3765366}}

\bibitem{unifsensitivity}
Hermon, J.
\newblock  On sensitivity of uniform mixing times.
\newblock {\em Ann. Inst. Henri Poincar\'e Probab. Stat.} \textbf{54} (2018), no.\ 1, 234--248. 
\newblock Available at: \href{https://projecteuclid.org/euclid.aihp/1519030827}{\nolinkurl{projecteuclid/1519030827}}

\bibitem{HK}
Hermon, J., and Kozma, G.,
\newblock  Sensitivity of mixing times of Cayley graphs.
Arxiv preprint \href{https://arxiv.org/abs/2008.07517}{\textcolor{blue}{arXiv:2008.07517}}

\bibitem{L2}
Hermon, J., and Peres, Y.,
\newblock A characterization of $L_2$ mixing and hypercontractivity via hitting times and maximal inequalities.
\newblock \emph{Probab. Theory Related Fields} 170 (2018), no. 3-4, 769--800. \href{http://www.ams.org/mathscinet-getitem?mr=MR3773799}{\textcolor{blue}{MR3773799}}

\bibitem{hermonsen}
Hermon, J., and Peres, Y.,
\newblock On sensitivity of mixing times and cutoff.
\newblock \emph{Electron. J. Probab.} 23 (2018), Paper No. 25, 34 pp. \href{http://www.ams.org/mathscinet-getitem?mr=MR3779818}{\textcolor{blue}{MR3779818}}



\bibitem{Keilson}
Keilson, J., 
\newblock Markov chain models -rarity and exponentiality. Applied Mathematical Sciences, 28. Springer-Verlag, New York-Berlin, 1979. {\rm xiii}+184 pp. ISBN: 0-387-90405-0 \href{http://www.ams.org/mathscinet-getitem?mr=MR0528293}{\textcolor{blue}{MR0528293}}













 

\bibitem{LPW} Levin, D., and Peres, Y., (2017).
 \newblock {\em Markov chains and mixing times.} 
 American Mathematical Society, Providence, RI. With contributions by Elizabeth L.\ Wilmer and a chapter by James G.\ Propp and David B.\ Wilson. \href{http://www.ams.org/mathscinet-getitem?mr=MR3726904}{\textcolor{blue}{MR3726904}}

 
\bibitem{LP}
Lyons, R., and Peres, Y., \emph{Probability on trees and networks}. Cambridge Series in Statistical and Probabilistic Mathematics, 42. Cambridge University Press, New York, 2016. \href{http://www.ams.org/mathscinet-getitem?mr=MR3616205}{\textcolor{blue}{MR3616205}} 



\bibitem{CRW}
Oliveira, R., 
\newblock Mean field conditions for coalescing random
walks.
\newblock \emph{Ann. Probab.} 41 (2013), no. 5, 3420--3461. \href{http://www.ams.org/mathscinet-getitem?mr=MR3127887}{\textcolor{blue}{MR3127887}}

\bibitem{Olivehit}
Oliveira, R.,
\newblock  Mixing and hitting times for finite Markov chains.
\newblock \emph{Electron. J. Probab.} 17 (2012), no. 70, 12 pp. \href{http://www.ams.org/mathscinet-getitem?mr=MR2968677}{\textcolor{blue}{MR2968677}}



\bibitem{CRWrev}
Oliveira, R., 
\newblock On the coalescence time of reversible random walks.
\newblock \emph{Trans. Amer. Math. Soc.} \textbf{364} (2012), no. 4, 2109--2128. \href{http://www.ams.org/mathscinet-getitem?mr=2869200}{\textcolor{blue}{MR2869200}}


\bibitem{PS}
Peres, Y., and Sousi, P.,
\newblock Mixing times are hitting times of large sets.
\newblock \emph{J. Theoret. Probab.} 28 (2015), no. 2, 488--519. \href{http://www.ams.org/mathscinet-getitem?mr=MR0942764}{\textcolor{blue}{MR3370663}}

\bibitem{inter}
Peres, Y., Sauerwald, T., Sousi, P., and Stauffer, A., Intersection and mixing times for reversible chains. \emph{Electron. J. Probab.} 22 (2017), No. 12, 16 pp. \href{http://www.ams.org/mathscinet-getitem?mr=MR3613705}{\textcolor{blue}{MR3613705}} 
\bibitem{Salez}
Salez, J., Cutoff for non-negatively curved Markov chains.  (2021) \href{https://arxiv.org/abs/2102.05597}{\textcolor{blue}{arXiv 2102.05597}}


\end{thebibliography}
\end{document}